\documentclass{amsart}
\usepackage[utf8]{inputenc}
\usepackage[T1]{fontenc}
\usepackage{graphicx}
\usepackage{subfigure}
\usepackage{epsfig}
\usepackage{xcolor}
\usepackage{geometry}
\usepackage{float}
\usepackage{empheq,etoolbox}
\usepackage{hyperref}

\usepackage{amsthm}
\usepackage{amsmath,amssymb,exscale}

\usepackage{framed}

\usepackage{enumerate}

\usepackage{tcolorbox}

\newtheorem{teo}{Theorem}[section]

\newtheorem{defi}[teo]{Definition}
\newtheorem{obs}[teo]{Remark}
\newtheorem{ex}[teo]{Example}

\newtheorem{prop}[teo]{Proposition}

\newtheorem*{remark}{Remark}

\newcommand{\dem}{\noindent \textit{Proof.} $\,$}
\newcommand{\fim}{\hfill $\rule{2.0mm}{2.0mm}$ \\}    % marca��o de fim das demonstra��es

\newcommand{\h}{\mathcal{H}}
\newcommand{\G}{\mathcal{G}}
\newcommand{\F}{\mathcal{F}}

\newcommand{\Leg}{\operatorname{Leg}}
\newcommand{\Tang}{\operatorname{Tang}}
\newcommand{\Sing}{\operatorname{Sing}}

\title{Homogeneous Convex Foliations of degree 6}
\author{Carla Pracias and Maycol Falla Luza} 
\email{carla\_pracias@id.uff.br \and hfalla@id.uff.br}
\date{}
\begin{document}
\begin{abstract}
In this paper, we study homogeneous convex foliations on the complex projective plane $\mathbb{P}^2$. A foliation is called convex if all of its leaves, except straight lines, have no inflection points, and such foliations form a Zariski closed subset in the space of degree $d$ foliations on $\mathbb{P}^2$. Using projective duality, every foliation can be associated with a $d$-web on the dual plane via its Legendre transform, and it is known that the Legendre transform of a homogeneous convex foliation is flat. Our first main result provides a classification of homogeneous convex foliations admitting exactly three radial singularities on the line at infinity. As a second result, we complete the classification of convex homogeneous foliations of degree $6$, extending previous classifications in degrees $4$ and $5$.
\end{abstract}

%==========================================================================
% Início do texto do artigo
%==========================================================================

\maketitle

\begin{quote}
\end{quote}

\begin{quote}
\end{quote}

\section{Introduction}

A \textbf{ foliation} on $\mathbb{P}^2$ is defined in affine coordinates by a polynomial $1$-form $\omega=a(x,y)dx+b(x,y)dy$ on $\mathbb{C}^2$ with isolated zeros. Equivalently the foliation can be defined by the vector field $X=b(x,y)\frac{\partial}{\partial x} - a(x,y) \frac{\partial}{\partial y} $. We say that a foliation is \textbf{convex} if all its leaves, except straight lines, have no inflection points. The set of convex foliations forms a Zariski closed subset of the space of degree $d$ foliations on $\mathbb{P}^2$. It can be seen that any foliation is convex if $d=0$ or $d=1$.  However, for $d \geq 2$, this subset is proper. In this note, we focus on homogeneous convex foliations. 

We say that $\h$ is a \textbf{homogeneous foliation} of degree $d$ on $\mathbb{P}^2$ if there exists a choice of affine coordinates $(x,y)$ where $\h$ is defined by a homogeneous vector field $A(x,y)\frac{\partial}{\partial x}+B(x,y)\frac{\partial}{\partial y}$ with $A,B\in \mathbb{C}[x,y]_d$ and $gcd(A,B)=1$. 

Using projective duality between lines and points, one can associate to any foliation $\F$ of degree $d$ on $\mathbb{P}^2$ a $d$--web on $\check{\mathbb{P}}^2$, called the Legendre transform and denoted by $\Leg \F$, as follows. Fix a generic line $l$ on $\mathbb{P}^2$ and consider the tangency locus $\Tang(\F,l)=\{p_1,\cdots,p_d\}\subset\mathbb{P}^2$ between $\F$ and $l$. We can think $l$ as a point of $\Check{\mathbb{P}}^2$ and the dual $\Check{p}_i$ as straight lines on $\Check{\mathbb{P}}^2$ passing through the point $l$. Then the set of tangent lines of $\Leg\F$ at $l$ is just $T_{\Check{l}}\Leg\F=\cup_{i=1}^d\Check{p}_i$. More precisely, let $(p,q)$ be the affine chart of $\Check{\mathbb{P}}^2$ correspond to the line $\{y=px+q\}\subset\mathbb{P}^2$. Then $\Leg \F$ is given by the implicit differential equation $\Check{F}(p,q;x):=a(x,px+q)-pb(x,px+q)=0$, with $x=-\frac{dq}{dp}$.

A distinguished class of webs on $\mathbb{P}^2$ consists of those with vanishing curvature, known as \textbf{flat webs}. Flat webs form a Zariski closed space that contains every web of maximal rank, see  for example \cite{H}. These webs have been extensively studied by many authors, and it is shown in \cite[Corollary 3.4]{tissus} that if $\h$ is a homogeneous convex foliation on $\mathbb{P}^2$ then $\Leg(\h)$ is flat.

Our first main result, Theorem \ref{Grau do tipo 3}, provides a classification of homogeneous convex foliations that have three radial singularities on the line at infinity. Radial singularities play a key role, as they correspond to fixed critical points of the Gauss map, as we shall see.

Homogeneous convex foliations of degrees $d = 4$ and $d = 5$ have been classified in \cite{bedrouni2020convex} and \cite{bedrouni2021convex}, respectively. Our second result, Theorem \ref{grau-6}, classifies convex homogeneous foliations of degree $d=6$.

\section{Basic definitions}
In this section, we review the main concepts needed for the rest of the article. For further details, the reader is referred to \cite{tissus}, \cite{JD}, and the references therein.

\subsection{Inflection divisor of a Foliation}

The \textbf{inflection divisor} of a foliation $\h$ of degree $d$ given by the homogeneous vector field $Z$ on $\mathbb{C}^3$, denoted by $I_\h$, is the curve formed, beside the singularities, by the inflection points of the leaves of $\h$. This divisor is defined by the vanishing of the determinant 
$$\det \left (\begin{matrix}
    x && y && z\\
    Z(x) && Z(y) && Z(z)\\
    Z^2(x) && Z^2(y) && Z^2(z)
\end{matrix} \right )$$
and has degree exactly $3d$. We consider the case when $\h$ is homogeneous, that is, we assume that in the affine chart $(x,y)$ of $\mathbb{P}^2$, $X=A\partial_x + B\partial_y$ is a homogeneous vector field. In this case the inflection divisor $I_\h$ decompose as $I_\h = C_{\h} + D_{\h} + L_{\infty}$ where
$$C_{\h} = xA + yB \hspace{1.5cm} \text{and} \hspace{1.5cm} D_{\h} = A_{x}B_{y} - A_{y}B_{x}$$
with $C_{\h} \in \mathbb{C}[x,y]_{d+1}$ and $D_{\h} \in \mathbb{C}[x,y]_{2d-2}$. On the other hand, $I_\h$ is decomposed by $I_\h=I_\h^{inv}+I_\h^{tr}$
where $I_\h^{inv}$ and $I_\h^{tr}$ are invariant and transverse part of $I_\h$, respectively. We have that:
\begin{itemize}
    \item $|I^{inv}_\h|$ consists of the lines of the tangent cone $C_\h$ and the infinity line $L_\infty$;
    \item $I_\h^{tr}=\sum_{i=1}^n (\rho_i-1)T_i$ for some number $n\le \deg D_\h$ of the lines $T_i$ passing through the origin $O$ of $\mathbb{C}^2$, $\rho_i-1$ being the inflection order of $T_i$.
\end{itemize}

 A homogeneous foliation $\h$ is \textbf{convex} if its inflection divisor is completely invariant by $\h$.

\subsection{The Gauss Map}

The \textbf{Gauss map} of a convex homogeneous foliation $\h$ is defined by rational map $\G_\h:\mathbb{P}^2\dashrightarrow\Check{\mathbb{P}}^2$ with $\G_\h(p)=T_p\h$, which is well defined outside the singular set $\Sing(\h)$. 

Since $\h$ is homogeneous we can also associate the rational map $\underline{\G}_\h:\mathbb{P}^1\to\mathbb{P}^1$ defined by
$$\underline{\G}_\h([y:x])=[-A(x,y):B(x,y)],$$

and this allows to determine completely the divisor $I_{\h}^{tr}$ and the set $\Sigma_{\h}^{rad}$ of radial singularities:
\begin{itemize}
    \item $\Sigma_{\h}^{rad}$ is formed of $[b:a:0]\in L_\infty$ such as $[a:b]\in\mathbb{P}^1$ is a fixed critical point of $\underline{\G}_\h$;
    \item $I_\h^{tr}=\sum_i(\rho_i-1) T_i$ where $T_i=(b_iy-a_ix=0)$ and $[a_i:b_i]\in\mathbb{P}^1$ is a non-fixed critical point of $\underline{\G}_\h$ of multiplicity $\rho_i-1$.
\end{itemize}

\subsection{Type of Homogeneous Foliation}

For a homogeneous foliation, not necessarily convex, we define the following.

\begin{defi}[\cite{tissus}, Definition 2.3]
Let $\h$ be a homogeneous foliation of degree $d$ on $\mathbb{P}^2_{\mathbb{C}}$
having a certain number $m\le 2d-2$ of
radial singularities $s_i$ of order $\tau_i -1$, $2\le \tau_i \le d$ for $i=1,2,\cdots,m$. The support of the divisor $I^{tr}_\h$ is constituted
of a certain number $n \le 2d-2$ of transverse inflection lines $T_j$ of order $\rho_j - 1$, $2\le \rho_j \le d$ for $j=1,2,\cdots,n$. We define the \textbf{type of the foliation} $\h$ by polynomial
$$\mathcal{T}_\h=\sum_{i=1}^{m}R_{\tau_i-1}+\sum_{j=1}^{n}T_{\rho_j-1}=\sum_{k=1}^{d-1}(r_k\cdot R_k+t_k\cdot T_k)$$
in $\mathbb{Z}[R_1,\cdots,R_{d-1},T_1,\cdots,T_{d-1}]$.
\end{defi}

Another useful definition for us is the \textbf{degree of type} defined by $\deg\mathcal{T}_\h=\sum\limits_{k=1}^{d-1}(r_k+t_k)$. Thus in our situation the type of a convex homogeneous foliation is only given by its invariant part, i.e., $\mathcal{T}_\h=\sum\limits_{k=1}^{d-1}(r_k\cdot R_k)$  and $\deg\mathcal{T}_\h=r_1+\hdots+r_k$.\\

\section{Homogeneous Convex Foliations with degree of type 3}

We consider on this section homogeneous convex foliations of degree $d\geq 3$ such that $\deg \mathcal{T}_{\h}=3$. For each $d$, we determine the number $N(d)$ of homogeneous convex foliations such that $\deg \mathcal{T}_{\h}=3$ and show that each type corresponds, up to conjugacy, to a unique foliation. The number of such foliations will be
$$N(d) = \left \{ \begin{matrix}
   \left \lfloor \frac{d+1}{3}\right \rfloor + \displaystyle \sum_{l=1}^{\lfloor d/3 \rfloor} \left \lfloor \frac{d-3l+2}{2} \right \rfloor -1 &,& \text{if}\quad d \not\equiv 1 \pmod{3}\\
   \left \lfloor \frac{d+1}{3}\right \rfloor + \displaystyle \sum_{l=1}^{\lfloor d/3 \rfloor} \left \lfloor \frac{d-3l+2}{2} \right \rfloor &,& \text{if} \quad d \equiv 1 \pmod{3}\end{matrix} \right .$$

Our main result here is to classify these foliations.

\begin{teo}\label{Grau do tipo 3}
Up to linear conjugacy, there are $N(d)$ homogeneous convex foliations of degree $d$ such that $\deg \mathcal{T}_{\h}=3$. They are described as follow:

\begin{enumerate}[i.]
    \item For each type $\mathcal{T}_{\h}=1\cdot R_{2l-2}+2\cdot R_{d-l}$, with $2\leq l \leq \dfrac{d+1}{3}$; or $\mathcal{T}_{\h}=3\cdot R_{d-l}$, when $d=3l-2$, there are unique homogeneous polynomials $A^{d}_{l}$, $B^{d}_{l}$ such that the foliation is given by $\omega^{d}_{l}=A^{d}_{l}dx + B^{d}_{l}dy$.

    \item For each type $\mathcal{T}_{\h}=1\cdot R_{\nu} + 1\cdot R_{d-\nu +l -2}+1\cdot R_{d-l}$, with $1\le l \le \frac{d}{3}$ and $2l-1 \le \nu \le \frac{d+l-2}{2}$, there are unique homogeneous polynomials $C^{d}_{l,\nu}$,  $D^{d}_{l,\nu}$ such that the foliation is given by $\omega^{d}_{l,\nu}=C^{d}_{l,\nu}dx + D^{d}_{l,\nu}dy$.
\end{enumerate}    
\end{teo}

Now we describe the elements of the theorem . For each integer $2\leq l \leq \dfrac{d+1}{3}$ we introduce the polynomials

\begin{align*}
A^{d}_{l}(x,y) &= 
a_{d-l+1}x^{l-1}y^{d-l+1}
+ a_{d-l+2}x^{l-2}y^{d-l+2}
+ \cdots + a_d y^d, \\[6pt]
B^{d}_{l}(x,y) &= 
b_0x^d + b_1x^{d-1}y + \cdots + b_{l-2}x^{d-l+2}y^{l-2} + x^{d-l+1}y^{l-1},
\end{align*}
where the coefficients 

\[
X =
\begin{pmatrix}
b_0 \\[4pt]
\vdots \\[4pt]
b_{\,l-2}
\end{pmatrix},
\qquad 
Y =
\begin{pmatrix}
a_{d-l+1} \\[4pt]
\vdots \\[4pt]
a_d
\end{pmatrix},
\]
\medskip
are determined by the following linear system
\[
T \cdot \begin{bmatrix} X \\[4pt] Y \end{bmatrix}
= \begin{bmatrix} C \\[4pt] D \end{bmatrix},
\qquad \text{where the matrix}\,\,\,
T=\begin{bmatrix} M & N \\[4pt] 0 & P \end{bmatrix}
\]
have blocks defined as follows:
\begin{itemize}
\item The matrix $M\in \mathbb{R}^{(l-1)\times(l-1)}$ is
\[
M=
\begin{pmatrix}
1 & 1 & 1 & \cdots & 1 \\[6pt]
0 & 1 & \tfrac{2!}{1!} & \cdots & \tfrac{(l-2)!}{(l-3)!}\\[8pt]
0 & 0 & 2! & \cdots & \tfrac{(l-2)!}{(l-4)!}\\[8pt]
\vdots & \vdots & \vdots & \ddots & \vdots \\[6pt]
0 & 0 & 0 & \cdots & (l-2)!
\end{pmatrix};
\]
\item The matrix $N\in \mathbb{R}^{(l-1)\times l}$ has entries
\[
N_{k+1,t}=\frac{(d-l+t)!}{(d-l+t-k)!},
\quad k=0,\dots,l-2,\;\; t=1,\dots,l;
\]
\item The matrix $P\in \mathbb{R}^{l\times l}$ has entries
\[
P_{i,t+1}=\frac{(d-l+1+t)!}{(d-2l+3+t-i)!},
\quad i=1,\dots,l,\;\; t=0,\dots,l-1.
\]
\end{itemize}

Finally, the vectors $C$ and $D$ are given by
\[
C=
\begin{pmatrix}
-1 \\[4pt]
-\tfrac{(l-1)!}{(l-2)!} \\[4pt]
\vdots \\[4pt]
-\tfrac{(l-1)!}{1!}
\end{pmatrix},
\qquad
D=
\begin{pmatrix}
-(l-1)! \\[4pt]
0 \\[4pt]
\vdots \\[4pt]
0
\end{pmatrix}.
\]

\noindent In the same way, for integers $1\le l \le \tfrac{d}{3}$ and $2l-1 \le \nu \le \tfrac{d+l-2}{2}$ we define
\begin{align*}
C^{d}_{l,\nu}(x,y) &= 
a_{d-\nu+l-1}x^{\nu-l+1}y^{d-\nu+l-1}
+ a_{d-\nu+l}x^{\nu-l}y^{d-\nu+l}
+ \cdots + a_d y^d, \\[6pt]
D^{d}_{l,\nu}(x,y) &= 
b_0x^d + b_1x^{d-1}y + \cdots + x^{d-l+1}y^{l-1}.
\end{align*}

\noindent where the coefficients 
\[
X=\begin{pmatrix} b_0 \\ \vdots \\ b_{l-2}\end{pmatrix},
\qquad 
Y=\begin{pmatrix} a_{d-\nu+l-1}\\ \vdots \\ a_d\end{pmatrix}
\]
are determined by a similar linear system
$
T\cdot \begin{bmatrix} X \\[4pt] Y \end{bmatrix}
=\begin{bmatrix} C \\[4pt] D \end{bmatrix},
$
with the follow modification
\begin{itemize}
\item The matrix $N\in \mathbb{R}^{(l-1)\times(\nu-l+2)}$ has entries
\[
N_{k+1,t+1}=\frac{(d-\nu+l-1+t)!}{(d-\nu+l-1+t-k)!},
\quad k=0,\dots,l-2,\;\; t=0,\dots,\nu-l+1;
\]

\item The matrix $P\in \mathbb{R}^{(\nu-l+2)\times(\nu-l+2)}$ has entries
\[
P_{i,t+1}=\frac{(d-\nu+l-1+t)!}{(d-\nu+t-i+1)!},
\quad i=1,\dots,\nu-l+2,\;\; t=0,\dots,\nu-l+1.
\]
\end{itemize}

%$$ 
%C^{d}_{l,\nu}=\frac{(d-\nu)}{(\nu - l+1)}\left[p_0x^{\nu-l+1}y^{d-\nu+l-1} +\sum_{i=1}^{\nu-l+1}(-1)^{i}p_ix^{\nu-l-i+1}y^{d-\nu+l+i-1}\right]
%$$
%with $p_0=\frac{(d-\nu+l)}{(d-\nu)}$, $p_1=-\binom{\nu-l+1}{1}$ and $p_{i}=\binom{\nu -l+1}{i}\prod_{j=2}^{i}\frac{(d-\nu+j-1)}{(d-\nu+l+j-1)}$ if $i\geq 2$, when $\nu -l$ is even, and

%$$
 %   C^{d}_{l,\nu}=\frac{(d-\nu)}{(\nu-l)(\nu-l+1)}\left[p_0x^{\nu-l+1}y^{d-\nu+l-1} + \sum_{i=1}^{\nu-l+1}(-1)^{i}p_ix^{\nu-l-i+1}y^{d-\nu+l+i-1}\right]
%$$  
%with $p_0=-\frac{(d-\nu+l+1)(d-\nu+l)}{(d-\nu)}$, $p_1 = -\binom{\nu -l +1}{1}(d-\nu+l+1)$, $p_2=-\binom{\nu -l +1}{2}(d-\nu+1)$ and $p_i=\binom{\nu -l+1}{i}\prod_{j=3}^{i}\frac{(d-\nu+j-1)}{(\nu-l+1)(d-\nu+l+j-1)}$, $i\geq 3$, when $\nu -l$ is odd. Finally we set the polynomial
    
%$$D^{d}_{l,\nu}= x^{d-l+1}y^{l-1} + \sum_{i=2}^{l}(-1)^{i}\binom{l-1}{i-1}\prod_{j=2}^{i}\frac{(d-\nu+j-1)}{(d-l+j)}x^{d-l+i}y^{l-i}.$$

\begin{ex}
Take $d=5$ and consider the foliations given by the following one-forms.
\begin{eqnarray*}
\h_{2}^5 : & \omega_{2}^5 = A_2^5 dx + B_2^5 dy =& \left(-xy^4 + \frac{3}{5}y^5 \right)dx + \left(x^4y + \frac{3}{5}x^5 \right)dy,\\
\h_{1,1}^5 : & \omega_{1,1}^5 = C_{1,1}^5 dx + D_{1,1}^5 dy =& \left(5xy^4 +4y^5 \right)dx + x^5 dy,\\
\h_{1,2}^5 : & \omega_{1,2}^5 = C_{1,2}^5 dx + D_{1,2}^5 dy =& \left(-10x^2y^3 +15xy^4 -6y^5 \right)dx + x^5 dy.
\end{eqnarray*}
A direct computation shows that $\mathcal{T}_{\h_{2}^5} = 1\cdot R_2 + 2\cdot R_3$, $\mathcal{T}_{\h_{1,1}^5} = 1\cdot R_1 + 1\cdot R_3 + 1\cdot R_4$ and $\mathcal{T}_{\h_{1,2}^5} = 2\cdot R_2 +1 \cdot R_4$, moreover, we can see that, up to linear conjugacy, these are the only homogeneous convex foliations of degree $5$ with degree of type $3$. In fact this is part of our main theorem.
\end{ex}

\noindent \textit{Proof of Theorem 3.1.} By hypothesis, $\h$ satisfies the following equations:
$$\left \{ \begin{matrix}
    r_1 + 2r_2 + \cdots + (d-1)r_{d-1}=2d-2 &,& r_{d-1} \le 1\\
    r_1 + r_2 + \cdots + r_{d-1}=3 &&
\end{matrix}
\right .$$

Given $1 \le l \le d-1$, note that $0\le r_{d-l} \le 3$. Assume that $r_{d-l}\neq 0$ and $r_j=0$ for $j>d-l$. We have three situations:

\begin{enumerate}[1.]
    \item if $r_{d-l}=3$ then we will have a foliation of type $\mathcal{T}_{\h}=3\cdot R_{d-l}$, but this only occurs if $3(d-l)=2d-2$, that is $d=3l-2$.

    \item if $r_{d-l}=2$, the foliation is of type $\mathcal{T}_{\h}=1\cdot R_{2l-2} +2\cdot R_{d-l}$ since $2d-2 -2(d-l)=2l-2$. Moreover, as $2l-2 \le d-l-1$ then $l \le \frac{d+1}{3}$. Therefore, we have a total of $\lfloor \frac{d+1}{3} \rfloor -1$ foliations of this type, where $2 \le l \le \frac{d+1}{3}$.

    \item if $r_{d-l}=1$, We identify the foliation by its type $\mathcal{T}_{\h}=1\cdot R_{\nu} + 1\cdot R_{d-\nu+l-2} + 1 \cdot R_{d-l}$, that is $2d-2-(d-l)-\nu=d-\nu+l-2$, where $\nu$ satisfies $2l-1\le \nu \le \frac{d+l-2}{2}$. From which we have $l \le \frac{d}{3}$. Note that for each $1 \le l \le \frac{d}{3}$ we have $K_d(l)=\lfloor \frac{d-3l+2}{2} \rfloor$ foliations and we obtain a total of $\sum_{l=1}^{\lfloor d/3 \rfloor} \left \lfloor \frac{d-3l+2}{2} \right \rfloor$ foliations of this type.
\end{enumerate}

To describe the foliations, let us consider the case where $\mathcal{T}_{\h}=1\cdot R_{2l-2} + 2\cdot R_{d-l}$:\\
Consider that the points $\infty = [1:0]$, $[0:1]$ and $[1:1]$ on $\mathbb{P}^{1}_{\mathbb{C}}$ are fixed critical points of $\mathcal{G}_{\h}$ with respective orders $d-l$, $d-l$ and $2l-2$. Thus, 

$B^{d}_l(1,z)=b_0 + b_1z + \cdots + z^{l-1}$,   $A^{d}_l(1,z)=a_{d-l+1}z^{d-l+1}+\cdots + a_dz^d$ and
 $h(z)=A^{d}_l(1,z) + B^{d}_l(1,z)=(z-1)^{2l-1}g_1(z)$.

Note that $l-1<d-l+1$ since $2l< 3l < d+1 <d+2$. Then, we have the following equality
        $$h(z)=b_0 +b_1z + \cdots b_{l-2}z^{l-2} + z^{l-1}+a_{d-l+1}z^{d-l+1}+a_{d-l+2}z^{d-l+2}+\cdots + a_dz^d = (z-1)^{2l-1}g(z)$$
which yields the system of $2l-1$ equations in $2l-1$ variables $h(1)=0, \, h'(1) =0,  \ldots \, , h^{(2l-2)}(1)=0$ that is precisely the system
\[
\begin{bmatrix}
M & N \\[6pt]
0 & P
\end{bmatrix}
\begin{bmatrix}
X \\[6pt]
Y
\end{bmatrix}
=
\begin{bmatrix}
C \\[6pt]
D
\end{bmatrix},
\]

We solve first the system $P\cdot Y = D$ and, since $M$ is triangular, this yields the full solution. Starting from the second row, we can factor out the common terms from the rows of the matrix, obtaining the equivalent system $\hat{P}\cdot Y = D$. It is easy to see that
$$
\det \hat{P}= d(d-1)^2\ldots(d-l+2)^{l-1}
\begin{vmatrix}
\displaystyle\prod_{m=0}^{\,l-2}\big(d-l+1-m\big) &
\displaystyle\prod_{m=0}^{\,l-3}\big(d-l+1-m\big) & \cdots &
d-l+1 & 1 \\[8pt]
\displaystyle\prod_{m=0}^{\,l-2}\big(d-l-m\big) &
\displaystyle\prod_{m=0}^{\,l-3}\big(d-l-m\big) & \cdots &
d-l & 1 \\[6pt]
\vdots & \vdots & \ddots & \vdots & \vdots \\[6pt]
\displaystyle\prod_{m=0}^{\,l-2}\big(d-2l+2-m\big) &
\displaystyle\prod_{m=0}^{\,l-3}\big(d-2l+2-m\big) & \cdots &
d-2l+2 & 1
\end{vmatrix}
$$
We claim that the last determinant is equal to $\displaystyle\prod_{j=1}^{l-1}j!$. In fact, denoting by $f_0(x) = 1$, $f_1(x) = x$, $f_2(x)= x(x-1), \ldots, f_{l-1}(x) = \displaystyle\prod_{j=0}^{l-2}(x-j)$ we have that the desired determinant is equal to
$$
\begin{vmatrix}
f_{l-1}(d-l+1) &
f_{l-2}(d-l+1) & \cdots &
f_1(d-l+1) & 1 \\[8pt]
f_{l-1}(d-l) &
f_{l-2}(d-l) & \cdots &
f_1(d-l) & 1 \\[6pt]
\vdots & \vdots & \ddots & \vdots & \vdots \\[6pt]
f_{l-1}(d-2l+2)&
f_{l-2}(d-2l+2) & \cdots &
f_1(d-2l+2) & 1
\end{vmatrix}.
$$

Alternatively, consider the polynomials $g_j(x)=x^j$, for $j=0, \ldots l-1$, and have two basis $\{f_j(x)\}$, $\{g_j(x)\}$ of the space of polynomials of degree at most $l-1$ and the transition matrix is triangular with ones on the diagonal. Hence our determinant reduces to a Vandermonde determinant:
$$
\begin{vmatrix}
g_{l-1}(d-l+1) &
g_{l-2}(d-l+1) & \cdots &
g_1(d-l+1) & 1 \\[8pt]
g_{l-1}(d-l) &
g_{l-2}(d-l) & \cdots &
g_1(d-l) & 1 \\[6pt]
\vdots & \vdots & \ddots & \vdots & \vdots \\[6pt]
g_{l-1}(d-2l+2)&
g_{l-2}(d-2l+2) & \cdots &
g_1(d-2l+2) & 1
\end{vmatrix}.
$$
which is equal to $ \displaystyle\prod_{1\leq i<j \leq l}((d-l+2-i) - (d-l+2-j))=\displaystyle\prod_{j=1}^{l-1}j!$. This proves that the system has a unique solution, as we wanted. For the case where $\mathcal{T}_{\h}=1\cdot R_{\nu} + 1\cdot R_{d-\nu +l -2}+1\cdot R_{d-l}$ we obtain a similar equation 
$$
h(z)=b_0 +b_1z + \cdots b_{l-2}z^{l-2} + z^{l-1}+a_{d-\nu+l-1}z^{d-\nu+l-1}+a_{d-\nu+l}z^{d-\nu+l}+\cdots + a_dz^d = (z-1)^{\nu+1}g(z)
$$
which we solve in an analogous way.
\fim

\begin{remark}
For polynomials $A^d_l$ and $B^{d}_{l}$ we can use similar computations in order to apply Cramer’s rule and find an explicit expression
$$A^{d}_{l}(x,y) \;=\; -\,x^{\,l-1}y^{\,d-l+1} 
+ \sum_{i=2}^{l}(-1)^{i}\binom{l-1}{i-1}
\prod_{j=2}^{i}\frac{d-2l+j}{\,d-l+j\,}\;
x^{\,l-i}y^{\,d-l+i},$$
$$B^{d}_{l}(x,y) \;=\; x^{\,d-l+1}y^{\,l-1} 
- \sum_{i=2}^{l}(-1)^{i}\binom{l-1}{i-1}
\prod_{j=2}^{i}\frac{d-2l+j}{\,d-l+j\,}\;
x^{\,d-l+i}y^{\,l-i}.$$

\noindent For example, we can compute 
\[
a_{d-l+2} = \frac{\det(\hat{P}_1)}{\det(\hat{P})},
\]
where $\hat{P}_1$ denotes the matrix obtained from $\hat{P}$ by replacing its first column with the vector $D$. Then 
$$
\det \hat{P}_1= -(l-1)!d(d-1)^2\ldots(d-l+2)^{l-1}
\begin{vmatrix}
\displaystyle\prod_{m=0}^{\,l-3}\big(d-l-m\big) & \cdots &
d-l & 1 \\[6pt]
\vdots & \vdots & \ddots & \vdots \\[6pt]
\displaystyle\prod_{m=0}^{\,l-3}\big(d-2l+2-m\big) & \cdots &
d-2l+2 & 1
\end{vmatrix}
$$
and the last determinant is computed as before, yielding the value $\displaystyle\prod_{j=1}^{l-2}j!$. Putting all these together we obtain 
$$
a_{d-l+2} = \frac{\det(\hat{P}_1)}{\det(\hat{P})}=\frac{-(l-1)!d(d-1)^2\ldots(d-l+2)^{l-1}\displaystyle\prod_{j=1}^{l-2}j!}{\det(\hat{P})}=-1.
$$

\noindent In general, the solutions are:

$a_{d-l+1}=-1, \hspace{0,2cm}
a_{d-l+2}= \binom{l-1}{1}\frac{(d-2l+2)}{(d-l+2)}, \hspace{0,2cm}
 \cdots, \hspace{0,2cm} a_{d}=(-1)^{l}\binom{l-1}{l-1}\frac{(d-l)\cdots(d-2l+2)}{d(d-1)\cdots(d-l+2)}$\\

$b_{l-2}= -\binom{l-1}{1}\frac{(d-2l+2)}{(d-l+2)}, \hspace{0,2cm} \hspace{0,2cm} \cdots, \hspace{0,2cm} b_{0}=(-1)^{l-1}\binom{l-1}{l-1}\frac{(d-l)\cdots(d-2l+2)}{d(d-1)\cdots(d-l+2)}.$\\

\end{remark}

Another nice observation is that we can provide a compact closed formula for the number $N(d)$.
\begin{prop}\label{N(d)}
For every integer $d \geq 1$, the quantity $N(d)$ is given by
\[
N(d) \;=\; \left\lfloor \frac{d-1}{2} \right\rfloor \;+\; \left\lfloor \frac{(d-1)^2+6}{12} \right\rfloor .
\]
\end{prop}

\dem
Define
\[
M(d):=\left\lfloor\frac{d-1}{2}\right\rfloor+\left\lfloor\frac{(d-1)^2+6}{12}\right\rfloor.
\]
We show that $M(d)=N(d)$ for all integers $d\geq 1$.  
The verification is carried out by distinguishing the residue $r\in\{0,1,2\}$ of $d \bmod 3$, and within each case analyzing the parity of $m=\lfloor d/3\rfloor$. We outline the general strategy and then give full details only in the case $r=0$ with $m$ even.

Write $d=3m+r$ with $r\in\{0,1,2\}$ and $m=\lfloor d/3\rfloor$. The definition of $N(d)$ involves the sum
\[
\sum_{i=1}^{m}\left\lfloor\frac{d-3i+2}{2}\right\rfloor
=\sum_{j=0}^{m-1}\left\lfloor\frac{3j+(r+2)}{2}\right\rfloor,
\]
where we substituted $j=m-i$. Setting $t:=r+2$, each summand can be written in the form
\[
\left\lfloor\frac{3j+t}{2}\right\rfloor=\frac{3j+t-\varepsilon_j}{2},
\qquad
\varepsilon_j\in\{0,1\},\ \varepsilon_j\equiv j+t\pmod 2.
\]
Summing over $j=0,\dots,m-1$ gives
\[
S=\frac{3m(m-1)}{4}+\frac{mt}{2}-\frac{E}{2},
\qquad
E:=\sum_{j=0}^{m-1}\varepsilon_j,
\]
and the value of $E$ depends only on the parity of $t$ and $m$.  
This leads to an explicit expression for $N(d)$ in terms of $m$, which can then be compared to the expression obtained for $M(d)$ by substituting $d=3m+r$ and expanding $(d-1)^2$. Equality is then checked case by case.

\medskip

\noindent\textbf{Case $r=0$ (i.e.\ $d=3m$), with $m$ even.}  
Since $r\neq 1$, the definition of $N(d)$ includes the adjustment $-1$, so
\[
N(3m)=m-1+S,\qquad 
S=\sum_{j=0}^{m-1}\left\lfloor\frac{3j+2}{2}\right\rfloor.
\]
For each $j$ we have
\[
\left\lfloor\frac{3j+2}{2}\right\rfloor=\frac{3j+2-(j\bmod 2)}{2}.
\]
Let $E=\sum_{j=0}^{m-1}(j\bmod 2)$, the number of odd indices among $0,\dots,m-1$. Then
\[
S=\frac{3m^2+m-2E}{4}.
\]

Now assume $m$ is even, say $m=2k$. In this case $E=\lfloor m/2\rfloor=k$, hence
\[
S=\frac{3(2k)^2+2k-2k}{4}=\frac{12k^2}{4}=3k^2.
\]
Therefore
\[
N(3m)=2k-1+3k^2=3k^2+2k-1.
\]

On the other hand, for $M(3m)$ we compute
\[
\left\lfloor\frac{3m-1}{2}\right\rfloor
=\left\lfloor\frac{6k-1}{2}\right\rfloor=3k-1,
\]
and
\[
\frac{(3m-1)^2+6}{12}=\frac{(6k-1)^2+6}{12}
=\frac{36k^2-12k+7}{12}=3k^2-k+\tfrac{7}{12},
\]
so
\[
\left\lfloor\frac{(3m-1)^2+6}{12}\right\rfloor=3k^2-k.
\]
Thus
\[
M(3m)=(3k-1)+(3k^2-k)=3k^2+2k-1,
\]
which coincides with $N(3m)$.

\medskip

The other cases ($r=0$ with $m$ odd, and $r=1,2$ with either parity of $m$) are entirely similar and yield the same conclusion. Therefore $M(d)=N(d)$ for all $d\geq 1$, proving the proposition.

\fim

\begin{obs}
The formula $N(d)$ we proved in Proposition \ref{N(d)} for the number of homogeneous convex foliations 
$\mathcal{H}$ of degree $d$ with $\deg T_{\mathcal{H}} = 3$, up to conjugacy, 
coincides with the number of connected unlabeled loopless multigraphs 
with three vertices and $d - 1$ edges. For further details, refer to sequence A253186 described in \textit{oeis.org}.
\end{obs}

\setcounter{secnumdepth}{2}
\section{Homogeneous Convex Foliations of Degree 6}

The aim of this section is to classify homogeneous convex foliations of degree six. The cases for degrees $d=4$ and $d=5$ have already been addressed in \cite{bedrouni2020convex} and \cite{bedrouni2021convex}, respectively. However, these classifications reveal that the complexity of the computations increases significantly with the degree. To handle the case $d=6$, we rely on primary decomposition of ideals, a technique implemented using SageMath. Some algebraic concepts and results are presented here, which will serve to establish the classification of homogeneous convex foliations of degree $6$.

\begin{defi}
Let $K$ be a field. If $I$, $J$ are ideals in $K[x_1,\hdots,x_n]$, then $I:J^{\infty}$ is the set
$$\{f\in K[x_1,\hdots,x_n]| \ \text{for all} \ g\in J, \ \text{there is} \ N\ge 0 \ \text{such that} \ fg^N \in I\}, $$
and is called the \textbf{saturation} of $I$ with respect to $J$.
\end{defi}

A result in the geometry context presented in \cite[Chapter 4]{ideals} says that $\overline{V(I)\setminus V(J)}\subseteq V(I:J^{\infty})$. Moreover, it is a straightforward computation that: if $J$ is a maximal ideal, then equality is verified.\\

The geometric interpretation of the primary decomposition is given by: Let $I = \underset{k}{\bigcap} I_k$ be a minimal primary decomposition, then $V(I) = \underset{k}{\bigcup} V(I_k)$.  If $I$ is a radical ideal, then each of the $I_k$ is a prime ideal minimal over $I$, and the primary decomposition simply expresses the algebraic set $V(I)$ as the union of the irreducible algebraic sets (algebraic varieties) $V(I_k)$. For further details on this computational method, see \cite[Section 3.8]{DP}.\\

\begin{teo}\label{grau-6} Let $\h$ be a homogeneous convex foliation of degree $6$ on $\mathbb{P}^2$, then $\h$ is conjugated to one of the following foliations:

\begin{enumerate}[1.]
    \item $\omega^6_2=y^6dx-x^6dy$;

    \item $\h^6_2$, $\h^6_{1,1}$, $\h^6_{1,2}$ and $\h^6_{2,3}$ as in Theorem \ref{Grau do tipo 3};

    \item for $\deg(\mathcal{T}_\h)=4$ the foliations are given by $\omega^6_{4,j}$, $j=1,\cdots,6$;

    \item for $\deg(\mathcal{T}_\h)=5$ the foliations are given by $\omega^6_{5,j,i}$,  $j=1,\cdots,6$, $i\in\{1,\cdots,6\}$;

    \item for $\deg(\mathcal{T}_\h)=6$ the foliations are given by $\omega^6_{6,j,i}$, $j=1,\cdots,5$,  $i\in\{1,\cdots,6\}$.
\end{enumerate}

where the forms are described in the Section \ref{solutions}.
\end{teo}

\dem
Let $\h$ be a homogeneous convex foliation of degree $d=6$ defined by the form $\omega=Adx + Bdy$ where $A,B \in \mathbb{C}[x,y]_d$. So it shall satisfies the equations
\begin{eqnarray*}
    \left \{ \begin{matrix}
        r_1 + 2r_2 + 3r_3 + 4r_4 + 5r_5 = 10 \\
        r_1 + r_2 + r_3 + r_4 + r_5 = \deg(\mathcal{T}_\h)
    \end{matrix} \right.
\end{eqnarray*}

with $\deg(\mathcal{T}_\h) \in \{2,3,4,5,6\}$. The case $\deg (\mathcal{T}_\h)=2$ was solved in \cite[Section 4]{tissus} so we have the following cases to consider

{\small
\begin{center}
\begin{tabular}{|c|c|c|c|c|}
    \hline
    4.1. $\deg(\mathcal{T}_\h)=3$ & 4.2. $\deg(\mathcal{T}_\h)= 4$ & 4.3. $\deg(\mathcal{T}_\h)= 5$ & 4.4. $\deg(\mathcal{T}_\h)= 6$\\
    \hline
    $\begin{matrix}
        1\cdot R_1 + 1\cdot R_4 + 1\cdot R_5 \\
        1\cdot R_2 + 1\cdot R_3 + 1\cdot R_5 \\
        1\cdot R_2 + 2\cdot R_4 \\
        2\cdot R_3 + 1\cdot R_4
    \end{matrix}$ & $\begin{matrix}
        2\cdot R_1 + 1\cdot R_3 + 1\cdot R_5 \\
        1\cdot R_1 + 2\cdot R_2 + 1\cdot R_5 \\
        2\cdot R_1 + 2\cdot R_4 \\
        3\cdot R_2 + 1\cdot R_4 \\
        1\cdot R_1 + 1\cdot R_2 + 1\cdot R_3 + 1\cdot R_4 \\
        2\cdot R_2 + 2\cdot R_3
    \end{matrix}$ &
    $\begin{matrix}
        3\cdot R_1 + 1\cdot R_2 + 1\cdot R_5 \\
        3\cdot R_1 + 2\cdot R_3 + 1\cdot R_4 \\
        2\cdot R_1 + 2\cdot R_2 + 1\cdot R_4 \\
        2\cdot R_1 + 1\cdot R_2 + 2\cdot R_3 \\
        1\cdot R_1 + 3\cdot R_2 + 1\cdot R_3 \\
        5\cdot R_2
    \end{matrix} $ & $\begin{matrix}
        5\cdot R_1 + 1\cdot R_5 \\
        4\cdot R_1 + 1\cdot R_2 + 1\cdot R_4 \\
        4\cdot R_1 + 2\cdot R_3 \\
        3\cdot R_1 + 2\cdot R_2 + 1\cdot R_3 \\
        2\cdot R_1 + 4\cdot R_2
    \end{matrix} $\\
    \hline
\end{tabular}
\end{center}}

\vspace{0.5cm}

\subsection{If  \texorpdfstring{$\deg(\mathcal{T}_{\h})=3$}{}}
In this case we have the foliations $\h^6_2$, $\h^6_{1,1}$, $\h^6_{1,2}$ and $\h^6_{2,3}$ with type $ 1\cdot R_2 + 2\cdot R_4 $, $1\cdot R_1 + 1\cdot R_4 + 1\cdot R_5$, $1\cdot R_2 + 1\cdot R_3 + 1\cdot R_5$, and $2\cdot R_3 + 1\cdot R_4$ respectively, given by Theorem \ref{Grau do tipo 3}.

\subsection{If  \texorpdfstring{$\deg(\mathcal{T}_\h)=4$}{}}

In the case $\mathcal{T}_\h=2\cdot R_1 + 1\cdot R_3 + 1\cdot R_5$, the map $\underline{\mathcal{G}}_\h$ has four fixed critical points (which we place at) $\infty = [1:0]$, $[0:1]$, $[1:1]$ and $[1:p]$ with $p\not\in \{0,1\}$ of orders 5, 3, 1 and 1, respectively. Thus, $A(1,z)=a_6z^6 + a_5z^5 + a_4z^4$ and $B(1,z)=1$ whose coefficients satisfy the following equations $Q_1: a_6 +a_5 +a_4+1=0$, $Q_2: 6a_6 + 5a_5 + 4a_4=0$, $Q_3: a_6p^6 + a_5p^5 + a_4p^4 + p=0$, $Q_4: 6a_6p^5 + 5a_5p^4 + 4a_4p^3=0$. \\
    We define the polynomial ring $R=\mathbb{Q}[a_6,a_5,a_4,p]$ and the ideals $I=<Q_1, Q_2, Q_3, Q_4>$, $J_1=<p>$, $J_2=<p-1>$ in $R$ and
$$J=(a_5 + 2a_4 + 6, a_6 - a_4 - 5, 3a_4p - 2a_4 + 15p, 8a_4^2 + 180p^2 + 115a_4 + 150p + 405, 54p^3 + 9p^2 - 8a_4 - 21p - 81)$$
is defined by saturation $J=(I:J_1^{\infty}):J_2^\infty$. The primary decomposition of $J$ is given by the intersection of the ideals
{\small $$(2p^4 - p^3 - p^2 - p + 2)\cap(-\frac{27}{4}p^3 - \frac{9}{8}p^2 + a_4 + \frac{21}{8}p + \frac{81}{8})\cap(\frac{27}{2}p^3 + \frac{9}{4}p^2 + a_5 - \frac{21}{4}p - \frac{57}{4})\cap(-\frac{27}{4}p^3 - \frac{9}{8}p^2 + a_6 + \frac{21}{8}p + \frac{41}{8}).$$}
Through the primary decomposition of $J$, we find the desired solutions. Therefore, the foliation is defined by
    $$\omega^6_{4,1}=[(f_1(p)-81)x^2y^4-2(f_1(p)-57)xy^5+(f_1(p)-41)y^6]dx + 8x^6dy,$$
    where the function $f_1(p)=54p^3+9p^2-21p$ with $p$ such that $2p^4-p^3-p^2-p+2=0$.\\

In cases $\mathcal{T}_\h=1\cdot R_1 + 2\cdot R_2 + 1\cdot R_5$, $\mathcal{T}_\h=2\cdot R_1 + 2\cdot R_4$, $\mathcal{T}_\h=3\cdot R_2 + 1\cdot R_4$, $\mathcal{T}_\h=1\cdot R_1 + 1\cdot R_2 + 1\cdot R_3 + 1\cdot R_4$ and $\mathcal{T}_\h=2\cdot R_2 + 2\cdot R_3$ we use the same strategy and obtain, respectively, the following foliations:\\

{\small \begin{tabular}{lclr}
$\omega^6_{4,2}$&=&$[-(f_2(p)+16)x^3y^3+3(f_2(p)+11)x^2y^4-3(f_2(p)+8)xy^5+(f_2(p)+6)y^6]dx + x^6dy$ & (\ref{5.2})\\
$\omega^6_{4,3}$&=&$[-6(f_3(p)+24)xy^5+5(f_3(p)+25)y^6]dx + [f_3(p)x^6 + 30x^5y]dy$ & (\ref{5.3})\\
$\omega^6_{4,4}$&=&$[10f_4(p)x^3y^3+15g_4(p)x^2y^4+6h_4(p)xy^5 + 10i_4(p)y^6]dx + [j_4(p)x^6 + 30x^5y]dy$ & (\ref{5.4})\\
$\omega^6_{4,5}$&=&$[15(f_5(p)-41)x^2y^4-12(2f_5(p)-77)xy^5+10(f_5(p)-37)y^6]dx + [-(2f_5(p)-1)x^6$ & (\ref{5.5})\\
&& $+ 60x^5y]dy$ & \\
$\omega^6_{4,6}$&=&$[5f_6(p)x^2y^4+3g_6(p)xy^5+2h_6(p)y^6]dx + [i_6(p)x^6 + j_6(p)x^5y + 300x^4y^2]dy$ & (\ref{5.6})
\end{tabular}}    
    
where the functions are given in the Section \ref{solutions}.

\subsection{If \texorpdfstring{$\deg(\mathcal{T}_\h)=5$}{}}

For $\mathcal{T}_\h=3\cdot R_1 + 1\cdot R_2 + 1\cdot R_5$, we can assume that $\infty = [1:0]$, $[0:1]$ are fixed critical points of map $\underline{\mathcal{G}}_\h$ with orders 5 and 2, respectively and $[1:1]$, $[1:p]$, $[1:q]$, with $p,q \not\in\{0,1\}$ and $p\ne q$ are fixed critical points of order 1. Then, $A(1,z)=a_6z^6 + a_5z^5+a_4z^4+a_3z^3$ and $B(1,z)=1$ whose coefficients satisfy the equations $Q_1: a_6 +a_5+ a_4 + a_3+ 1=0$, $Q_2: 6a_6 + 5a_5 + 4a_4 + 3a_3=0$, $Q_3: a_6p^6 + a_5p^5 a_4p^4 + a_3p^3 + p=0$, $Q_4:  6a_6p^5 + 5a_5p^4 + 4a_4p^3 +3a_3p^2=0$, $Q_5: a_6q^6 + a_5q^5+ a_4q^4 + a_3q^3 + q=0$ and $Q_6: 6a_6q^5 + 5a_5q^4 + 4a_4q^3 +3a_3q^2=0$.

We define the polynomial ring $R=\mathbb{Q}[a_6,a_5,a_4,a_3,p,q]$ and the ideals $I=<Q_1, Q_2, Q_3, Q_4, Q_5,$ $ Q_6>$, $J_1=<p>$, $J_2=<p-1>$, $J_3=<q>$, $J_4=<q-1>$, $J_5=<p-q>$ in $R$ and $J=((((I:J_1^{\infty}):J_2^\infty):J_3^\infty):J_4^\infty):J_5^\infty$. The primary decomposition of $J$ gives us the foliation
$$\omega^6_{5,1}=[8f_1(q)x^3y^3+g_1(q)x^2y^4+2h_1(q)xy^5 +i_1(q)y^6]dx + 8x^6dy $$
where the function are given in the Section \ref{solutions} (\ref{5.7}).\\

In cases $\mathcal{T}_\h=3\cdot R_1 + 1\cdot R_3 + 1\cdot R_4$, $\mathcal{T}_\h=2\cdot R_1 + 2\cdot R_2 + 1\cdot R_4$, $\mathcal{T}_\h=2\cdot R_1 + 1\cdot R_2 + 2\cdot R_3$, $\mathcal{T}_\h=1\cdot R_1 + 3\cdot R_2 + 1\cdot R_3$, $\mathcal{T}_\h=5\cdot R_2$ we use the same strategy and obtain, respectively, the following foliations:\\

{\small \begin{tabular}{lclr}
$\omega^6_{5,2}$&=&$[f_2(q)x^2y^4+g_2(q)xy^5 +h_2(q)y^6]dx + [i_2(q)x^6 + x^5y]dy$ & (\ref{5.8})\\
$\omega^6_{5,3,i}$&=&$[f_3(q)x^3y^3+g_3(q)x^2y^4 +h_3(q)xy^5 +i_3(q)y^6]dx + [j_3(q)x^6 + x^5y]dy$ & (\ref{5.9}) \\
$\omega^6_{5,4,i}$&=&$[f_4(q)x^2y^4+g_4(q)xy^5 +h_4(q)y^6]dx + [i_4(q)x^6 + j_4(q)x^5y + x^4y^2]dy$ & (\ref{5.10})\\
$\omega^6_{5,5}$&=&$[f_5(q)x^3y^3+g_5(q)x^2y^4+h_5(q)xy^5 +i_5(q)y^6]dx + [j_5(q)x^6 + k_5(q)x^5y + x^4y^2]dy$ & (\ref{5.11})\\
$\omega^6_{5,6,i}$&=&$[f_6(q)x^3y^3+g_6(q)x^2y^4+h_6(q)xy^5 +i_6(q)y^6]dx + [j_6(q)x^6 + k_6(q)x^5y + l_6(q)x^4y^2 $ & (\ref{5.12})\\
&& $+x^3y^3]dy$ & 
\end{tabular}}

where the functions are given in the Section \ref{solutions}.

\subsection{If  \texorpdfstring{$\deg(\mathcal{T}_\h)=6$}{}}

Firstly, if $\mathcal{T}_\h=5\cdot R_1 + 1\cdot R_5$ the map $\underline{\mathcal{G}}_\h$ has one point fixed critical of order 5 and five of order 1, wich are $\infty = [1:0]$, $[0:1]$, $[1:1]$, $[1:p]$, $[1:q]$ and $[1:r]$, respectively, with $p,q,r \not\in\{0,1\}$ and $p, q,r$ are pairwise distinct. Thus, $A(1,z)=a_6z^6 + a_5z^5+a_4z^4+a_3z^3+a_2z^2$ and $B(1,z)=1$ whose coefficients satisfy the equations $Q_1: a_6 +a_5+ a_4 + a_3 + a_2 + 1=0$, $Q_2: 6a_6 + 5a_5 + 4a_4 + 3a_3 + 2a_2=0$, $Q_3: a_6p^6 + a_5p^5 + a_4p^4 + a_3p^3+ a_2p^2 + p=0$, $Q_4: 6a_6p^5 + 5a_5p^4 + 4a_4p^3 +3a_3p^2 +2a_2p=0$, $Q_5: a_6q^6 + a_5q^5 + a_4q^4 + a_3q^3+ a_2q^2 + q=0$, $Q_6: 6a_6q^5 + 5a_5q^4 + 4a_4q^3 +3a_3q^2 +2a_2q=0$, $Q_7: a_6r^6 + a_5r^5 + a_4r^4 + a_3r^3+ a_2r^2 + r=0$, $Q_8: 6a_6r^5 + 5a_5r^4 + 4a_4r^3 +3a_3r^2 +2a_2r=0$.
    
We define the polynomial ring $R=\mathbb{Q}[a_6,a_5,a_4,a_3,a_2,p,q,r]$, and the ideals $I=<Q_1, Q_2, Q_3, Q_4, $ $Q_5, Q_6, Q_7, Q_8>$, $J_1=<p>$, $J_2=<p-1>$, $J_3=<q>$, $J_4=<q-1>$, $J_5=<r>$, $J_6=<r-1>$, $J_7=<p-q>$ e $J_8=<q-r>$ in $R$. Given $J=((((I:J_1^{\infty}):J_2^\infty):\cdots):J_8^\infty$, we have that the solutions of equations is given by the primary decomposition of $J$, and the foliations follows:
$$\omega^6_{6,1,i}=[f_1(r)x^4y^2+g_1(r)x^3y^3+h_1(r)x^2y^4+i_1(r)xy^5 +j_1(r)y^6]dx + x^6dy$$
   where the functions are given in the Section \ref{solutions} (\ref{5.13}).\\

In cases $\mathcal{T}_\h=4\cdot R_1 + 1\cdot R_2 + 1\cdot R_4$, $\mathcal{T}_\h=3\cdot R_1 + 2\cdot R_2 + 1\cdot R_3$, $\mathcal{T}_\h=4\cdot R_1 + 2\cdot R_3$, $\mathcal{T}_\h=2\cdot R_1 + 4\cdot R_2$ we use the same strategy and obtain, respectively, the following foliations:\\

{\small \begin{tabular}{lclr}
$\omega^6_{6,2}$&=&$[f_2(r)x^3y^3+g_2(r)x^2y^4+h_2(r)xy^5 +i_2(r)y^6]dx + [j_2(r)x^6 +x^5y]dy$ & (\ref{5.14})\\
$\omega^6_{6,3}$&=&$[f_3(r)x^3y^3+g_3(r)x^2y^4+h_3(r)xy^5 +i_3(r)y^6]dx + [j_3(r)x^6 +k_3(r)x^5y +x^4y^2]dy$ & (\ref{5.15})\\
$\omega^6_{6,4,i}$&=&$[f_4(r)x^2y^4+g_4(r)xy^5 +h_4(r)y^6]dx + [i_4(r)x^6 +j_4(r)x^5y +x^4y^2]dy$ & (\ref{5.16})\\
$\omega^6_{6,5,i}$&=&$[f_5(r)x^3y^3+g_5(r)x^2y^4+h_5(r)xy^5 +i_5(r)y^6]dx + [j_5(r)x^6 +k_5(r)x^5y +l_5(r)x^4y^2$ & (\ref{5.17})\\
&& $+ x^3y^3]dy$ &
\end{tabular}}

where the functions are given in the Section \ref{solutions}.

%%%%%%%%%%%%%%%%%%%%%%%%%%%%%%%%%%%%%%%%%%%%%%%%%%%%%%%%%%%%%%%%%%%%%%%%%%%%%%%%%%%%%%%%%%%%%%%%%%%%%%%%%%%%%%%%%%%%%%%%%%%%%%%%%%%%%%%%%%%%%%%%%%%%%%%%%%%%%%%%%%%%%%%%%%%%%%%%%%%%

\section{Description of the Solutions}{\label{solutions}}

In what follows, we describe the forms given in Theorem \ref{Grau do tipo 3}.
\subsection{}\label{5.1} $\omega^6_{4,1}=[(f_1(p)-81)x^2y^4-2(f_1(p)-57)xy^5+(f_1(p)-41)y^6]dx + 8x^6dy$

$f_1(p)=54p^3+9p^2-21p$ with $p$ such that $2p^4-p^3-p^2-p+2=0$.

\subsection{}\label{5.2}$\omega^6_{4,2}=[-(f_2(p)+16)x^3y^3+3(f_2(p)+11)x^2y^4-3(f_2(p)+8)xy^5+(f_2(p)+6)y^6]dx + x^6dy$

$f_2(p)=16p^3-24p^2-4p$ with $p$ such that $(2p^2-2p+1)(p^2-p-1)=0$.

\subsection{}\label{5.3}$\omega^6_{4,3}=[-6(f_3(p)+24)xy^5+5(f_3(p)+25)y^6]dx + [f_3(p)x^6 + 30x^5y]dyy$

$f_3(p)=p^4+p^3+p^2-24p$ with $p$ such that $p^6+2p^5+3p^4-21p^3+3p^2+2p+1=0$.

\subsection{}\label{5.4}$\omega^6_{4,4}=[10f_4(p)x^3y^3+15g_4(p)x^2y^4+6h_4(p)xy^5 + 10i_4(p)y^6]dx + [j_4(p)x^6 + 30x^5y]dy$

\begin{tabular}{lcl}
$f_4(p)$ &=& $-2p^5+3p^4-28p^3+12p^2+39p-6$  \\
$g_4(p)$ &=& $2p^5-p^4+24p^3+19p^2-68p-5$ \\
$h_4(p)$ &=& $-2p^5-5p^4-12p^3-112p^2+155p+33$ \\ 
$i_4(p)$ &=& $2p^4-4p^3+31p^2-29p-9$ \\
$j_4(p)$ &=& $2p^5-5p^4+32p^3-43p^2+10p-3$
\end{tabular}

and $p$ such that $p^6 - 3p^5 + 18p^4 - 31p^3 + 18p^2 - 3p + 1=0$.

\subsection{}\label{5.5}$\omega^6_{4,5}=[15(f_5(p)-41)x^2y^4-12(2f_5(p)-77)xy^5+10(f_5(p)-37)y^6]dx + [-(2f_5(p)-1)x^6+ 60x^5y]dy$

$f_5(p)=6p^5+4p^4-17p^3+28p^2+42p$ with  $p$ is such that $6p^6-8p^5-7p^4+44p^3-37p^2+2p+1=0$

\subsection{}\label{5.6}$\omega^6_{4,6}=[5f_6(p)x^2y^4+3g_6(p)xy^5+2h_6(p)y^6]dx + [i_6(p)x^6 + j_6(p)x^5y + 300x^4y^2]dy$

{\small \begin{tabular}{lcl}
$f_6(p)$ &=& $13p^5-427p^4+689p^3-451p^2+760p-359$  \\
$g_6(p)$ &=& $-32p^5+1051p^4-1694p^3+1121p^2-1822p+795$ \\
$h_6(p)$ &=& $19p^5-624p^4+1005p^3-670p^2+1061p-446$ \\ 
$i_6(p)$ &=& $-8p^5+263p^4-430p^3+245p^2-622p+7$ \\
$j_6(p)$ &=& $p^5-33p^4+57p^3-13p^2+161p-5$ \\

\end{tabular}}

and $p$ such that $p^6-33p^5+58p^4-41p^3+58p^2-33p+1=0$.

\subsection{}\label{5.7}$\omega^6_{5,1}=[8f_1(q)x^3y^3+g_1(q)x^2y^4+2h_1(q)xy^5 +i_1(q)y^6]dx + 8x^6dy$

{\small \begin{tabular}{lcl}
$f_1(q)$ &=& $-32q^{11}+44q^{10}+82q^9-205q^8+306q^6-240q^5-101q^4+188q^3-24q^2-56q+7$ \\
$g_1(q)$ &=& $720q^{11}-900q^{10}-1950q^9-4377q^8+510q^7-6816q^6+4684q^5+2787-3954q^3$\\
&& $+192q^2+1272q-123$\\
$h_1(q)$ &=& $-336q^{11}+372q^{10}+966q^9-1917q^8-510q^7+3144q^6-1764q^5-1575q^4+1698q^3$\\
&&$+96q^2-600q+15$\\
$i_1(q)$ &=& $208q^{11}-196q^{10}-638q^9+1097q^8+510q^7-1920q^6+804q^5+1171q^4-946q^3$\\
&&$-192q^2+376q+29$
\end{tabular}}

and $q$ such that $4q^{12}-10q^{11}-3q^{10}+36q^9-32q^8-32q^7+75q^6-32q^5-32q^4+36q^3-3q^2-10q+4=0$.

\subsection{}\label{5.8}$\omega^6_{5,2}=[f_2(q)x^2y^4+g_2(q)xy^5 +h_2(q)y^6]dx + [i_2(q)x^6 + x^5y]dy$

{\small\begin{tabular}{lcl}
    $f_2(q)$&=& $-(6523/7500q^{17} + 49689/10000q^{16} - 546067/18000q^{15} - 4470013/36000q^{14}$\\
    &&$+ 24363877/36000q^{13} + 33240323/120000q^{12} - 141951433/30000q^{11} + 104569223/18000q^{10}$  \\
    && $ + 124404953/36000q^9 - 926170303/72000q^8 + 161735553/20000q^7 + 93004537/30000q^6$ \\
    && $- 106279669/18000q^5 + 128702149/72000q^4 + 18225181/36000q^3 - 19072591/60000q^2$ \\
    && $+ 147517/30000q + 17687/1200)$\\
    $g_2(q)$&=& $-(-65167/37500q^{17} - 372109/37500q^{16} + 13647193/225000q^{15}+ 3719519/15000q^{14} $\\
    &&$ - 13532217/10000q^{13} - 61969441/112500q^{12} + 709555867/75000q^{11} - 2616200957/225000q^{10}$ \\
    &&$  - 207149351/30000q^9 + 385999937/15000q^8 - 3641803607/225000q^7 - 38767209/6250q^6  $ \\
    && $+ 5316275411/450000q^5 - 53712013/15000q^4 - 3793739/3750q^3 + 71605709/112500q^2$\\
    &&  $ - 537109/37500q - 76491/3125)$\\
    $h_2(q)$&=& $-(521/600q^{17} + 111521/22500q^{16} - 2728841/90000q^{15} - 418034/3375q^{14}$\\
    &&$+ 73062199/108000q^{13} + 29631359/108000q^{12} - 425653001/90000q^{11} + 523502909/90000q^{10}$ \\
    && $  + 372638147/108000q^9 - 1389829319/108000*q^8  + 218618363/27000q^7 + 139598717/45000q^6$\\
    &&$ - 1063560677/180000q^5   + 193569929/108000q^4 + 27299687/54000q^3 - 17198729/54000q^2 $\\
    && $  + 5411/625q + 483217/45000)$\\
    $i_2(q)$&=& $-(-7/25000q^{17} - 1117/450000q^{16} + 83/25000q^{15} + 32951/540000q^{14}$\\
    &&$- 3679/67500q^{13} - 2849317/5400000q^{12} + 100649/225000q^{11} + 328397/225000q^{10}$ \\
    && $  - 72089/67500q^9 - 1147729/1080000q^8 + 5507329/2700000q^7 + 183823/450000q^6$ \\
    && $ - 763987/900000q^5+ 1033411/1080000q^4 - 76169/540000q^3 - 333971/2700000q^2 $\\
    && $ + 112211/150000q - 91/450000)$\\
    \end{tabular}}
    
and $q$ such that $12q^{18} + 72q^{17} - 398q^{16} - 1828q^{15} + 8815q^{14} + 6367q^{13} - 63493q^{12} +61822q^{11} + 65677q^{10}- 158515q^9 + 65677q^8 + 61822q^7 - 63493q^6 + 6367q^5 + 8815q^4 - 828q^3 - 398q^2 + 72q + 12=0$ 

\subsection{}\label{5.9}$\omega^6_{5,3}=[f_3(q)x^3y^3+g_3(q)x^2y^4 +h_3(q)xy^5 +i_3(q)y^6]dx + [j_3(q)x^6 + x^5y]dy$
    
{\small    \begin{tabular}{lcl}
    $f_3(q)$&=& $-(2073/5000q^{17} - 7073/625q^{16} + 244193/3000q^{15} - 3629657/9000q^{14}$\\
    &&$+ 51874721/36000q^{13} - 370187701/90000q^{12} + 1695170921/180000q^{11} - 289985123/18000q^{10}$ \\
    && $  + 1401146779/72000q^9 - 28547647/2250q^8 - 848904877/180000q^7 + 1232337767/60000q^6$ \\
    &&  $- 1647758459/72000q^5 + 107864197/7200q^4 - 39794417/6000q^3 + 335366183/180000q^2 $\\
    &&$- 9361357/40000q + 136/9)$\\
    $g_3(q)$&=& $-(-6219/10000q^{17} + 114657/4000q^{16} - 2156089/10000q^{15} + 2642947/2400q^{14}$\\
    &&$- 95975407/24000q^{13} + 2774737199/240000q^{12} - 161058161/6000q^{11} + 11164611539/240000q^{10}$ \\
    && $ - 2737992533/48000q^9 + 3675246629/96000q^8 +1476820757/120000q^7 - 190177623/3200q^6$\\
    &&$ + 8067959497/120000q^5 - 4254093481/96000q^4 + 315198673/16000q^3 - 83494991/15000q^2 $ \\
    && $ + 5528201/8000*q - 16947787/480000)$\\
    $h_3(q)$&=& $-(6219/25000q^{17} - 318657/12500q^{16} + 4961461/25000q^{15} - 15540499/15000q^{14}$\\
    &&$+ 45655493/12000q^{13} - 104010881/9375q^{12} + 7799140117/300000q^{11} - 3407417347/75000q^{10} $\\
    && $ + 1349705959/24000q^9 - 2304959573/60000q^8 - 3360568397/300000q^7 + 5811867849/100000q^6 $\\
    && $  - 39827299091/600000q^5 + 1318065263/30000q^4 - 97907711/5000q^3 + 1665741163/300000q^2 $\\
    && $- 136875827/200000*q + 9087787/300000)$\\
    $i_3(q)$&=& $-(77927/10000q^{16} - 77927/1250q^{15} + 5955421/18000q^{14} - 22050343/18000q^{13}$\\
    &&$+ 258797279/72000q^{12} - 1525992299/180000q^{11} + 5364909079/360000q^{10}$ \\
    && $ - 668422877/36000q^9 + 1848197221/144000q^8 + 15697897/4500q^7 - 2289765041/120000q^6$\\
    && $  + 7897126699/360000q^5 - 2096809541/144000q^4 + 31204201/4800q^3 - 66518749/36000q^2 $\\
    && $+ 27197939/120000q - 6547787/720000)$
    \end{tabular}
    
    \begin{tabular}{lcl}
    $j_3(q)$&=&$-(-2073/50000q^{17} + 35241/100000q^{16} - 285841/150000q^{15} + 144877/20000q^{14}$\\
    && $ -518269/24000q^{13} + 20751601/400000q^{12} - 28196437/300000q^{11} +
144931127/1200000q^{10}$\\
 	&& $- 857347/9600q^9 - 2349757/160000q^8 +
24554333/200000q^7 - 174910493/1200000q^6$\\
	&& $ + 9490711/100000q^5 -
20158133/480000q^4 + 1052221/80000q^3 - 462073/300000q^2$\\
	&& $ +221533/300000q - 1357/800000)$
    \end{tabular}}

and $q$ such that $i= \left \{\begin{matrix}
    1, && q^2-q+3=0\\
    2, && 72q^{16} - 576q^{15} + 2840q^{14} - 9800q^{13} + 26170q^{12} - 55828q^{11} \\
    && + 81274q^{10}- 59540q^9 - 14825q^8 + 79280q^7 - 85758q^6 \\
    &&+ 54454q^5 - 23135q^4 + 6030q^3 - 700q^2 + 42q - 1=0
\end{matrix} \right. $

\subsection{}\label{5.10}$\omega^6_{5,4,i}$=$[f_4(q)x^2y^4+g_4(q)xy^5 +h_4(q)y^6]dx + [i_4(q)x^6 + j_4(q)x^5y + x^4y^2]dy$

{\small  \begin{tabular}{lcl}
    $f_4(q)$&=& $-(363757/1350000q^{19} - 88348429/10800000q^{18} + 378994811/4050000q^{17}$\\
    &&$- 2336711699/10800000q^{16} - 5379624499/10800000q^{15} + 93868823629/16200000q^{14}$ \\
    && $ - 269880383269/10800000q^{13} + 90189476507/1350000q^{12} - 2149197428351/16200000q^{11}$ \\
    &&$ + 79572875673/400000q^{10} - 2265191830589/10800000q^9 + 4912544034293/32400000q^8 $\\
    &&$ - 142097994839/1800000q^7 + 3418823621/112500q^6 - 225325737547/32400000q^5$\\
    && $+ 198665033/5400000q^4 + 5531599061/10800000q^3 - 6438503657/32400000q^2$\\
    && $ + 46693859/2700000q + 51430871/10800000)$\\
    $g_4(q)$&=& $-(-246799/562500q^{19} + 9579109/720000q^{18} - 2050640551/13500000q^{17}$\\
    &&$+ 6246624503/18000000q^{16} + 14768588383/18000000q^{15} - 253770180247/27000000q^{14}$ \\
    && $ + 16147093241/400000q^{13} - 40315211201/375000q^{12} + 5748247613207/27000000q^{11}$ \\
    && $- 5726496629267/18000000q^{10} + 6005383761457/18000000q^9 - 2590606910401/10800000q^8$\\
    &&$ + 1120619384761/9000000q^7 - 107522990429/2250000q^6 + 586504621459/54000000q^5$\\
    && $- 593984339/9000000q^4 - 2869448933/3600000q^3 + 16584618221/54000000q^2$\\
    && $ - 15011311/562500q - 110202407/18000000)$\\
    $h_4(q)$&=& $-(69271/375000q^{19} - 1049509/187500q^{18} + 143627881/2250000q^{17}$\\
    && $- 108608689/750000q^{16} - 782413657/2250000q^{15} + 8883408701/2250000q^{14}$ \\
    &&$  - 14273087893/843750q^{13} + 151701790699/3375000q^{12} - 49986808123/562500q^{11} $\\
    &&$ + 223626811631/1687500q^{10} - 935047982717/6750000q^9 + 167510219119/1687500q^8$ \\
    &&$ - 173581350061/3375000q^7 + 66497106977/3375000q^6 - 15049120163/3375000q^5$\\
    && $+ 65886551/2250000q^4 + 2203911401/6750000q^3 - 845509547/6750000q^2$\\
    && $ + 8155181/750000q + 1261907/562500)$\\
    $i_4(q)$&=& $-(-51869/2250000q^{19} + 2480887/3600000q^{18} - 3243052/421875q^{17}$\\
    && $+ 269896837/18000000q^{16} + 879684757/18000000q^{15} - 12732080783/27000000q^{14} $\\
    && $ + 20767295069/10800000q^{13} - 32835151877/6750000q^{12} + 250169361553/27000000q^{11} $\\
    &&$ - 710655042679/54000000q^{10} + 684808369969/54000000q^9 - 29852198501/3600000q^8$ \\
    &&$ + 109768846937/27000000q^7 - 4958281799/3375000q^6 + 15460185061/54000000q^5$\\
    &&$ - 7209019/1000000q^4 - 41459761/2160000q^3 + 325572719/54000000q^2$\\
    &&$  - 2237227/4500000q + 249847/18000000)$\\
    $j_4(q)$&=& $-(12883/1687500q^{19} - 11645743/54000000q^{18} + 88004677/40500000q^{17}$\\
    && $- 46179917/54000000q^{16} - 1268769157/54000000q^{15} + 10359951709/81000000q^{14} $\\
    &&$  - 20814521383/54000000q^{13} + 2078994781/3375000q^{12} - 51163412813/81000000q^{11}$ \\
    && $- 2750419189/18000000q^{10} + 105383360341/54000000q^9 - 441248618329/162000000q^8 $\\
    &&$ + 5388191317/3000000q^7 - 1879365443/2250000q^6 + 43092035999/162000000q^5$\\
    &&$ + 192632753/27000000q^4 - 1211058493/54000000q^3 + 1754174593/162000000q^2$\\
    && $ - 829963/843750q + 5560253/54000000)$
    \end{tabular}}

and $q$ such that $i= \left \{\begin{matrix}
    1, && q^4 - 2q^3 + 7q^2 - 2q + 1=0\\
    2, && 3q^{16} - 87q^{15} + 905q^{14} - 645q^{13} - 11870q^{12} + 50818q^{11} \\
    &&- 136452q^{10} + 266320q^9 - 337985q^8 + 266320q^7 - 136452q^6 \\
    &&+ 50818q^5 - 11870q^4 - 645q^3 + 905q^2 - 87q + 3=0
\end{matrix} \right. $

\subsection{}\label{5.11}$\omega^6_{5,5}=[f_5(q)x^3y^3+g_5(q)x^2y^4+h_5(q)xy^5 +i_5(q)y^6]dx + [j_5(q)x^6 + k_5(q)x^5y + x^4y^2]dy$

{\small    \begin{tabular}{lcl}
    $f_5(q)$&=& $-(534056809401/3012275000q^{17} - 3264135724403/3012275000q^{16}$ \\
    &&$+ 4857291489431/753068750q^{15} - 6716963422/307375q^{14} + 999079511509/15061375q^{13}$\\
    &&$- 100067168296801/753068750q^{12} + 60542349319849/376534375q^{11}$\\
    &&$+ 190625466597421/4518412500q^{10} - 271474648130507/361473000q^9$ \\
    && $ + 1122220235939109/602455000q^8 - 2114523954168151/753068750q^7 $\\
    && $+ 12354053474584963/4518412500q^6 - 2214856089985783/1506137500q^5$ \\
    && $+ 123905697365281/451841250q^4 + 11613786594307/225920625q^3$\\
    &&$ - 21143890937117/1506137500q^2 + 7491217119661/9036825000q $\\
    &&$- 84921525293/9036825000)$\\
    $g_5(q)$&=& $-(-194440858887/172130000q^{17} + 3858579575169/430325000q^{16}$\\
    &&$-3097235585331/53790625q^{15} + 3015196417659/12295000q^{14} $ \\
    && $- 75393011669463/86065000q^{13} + 21255032948931/8606500q^{12} $ \\
    && $- 2417326884708161/430325000q^{11} + 2199355061717921/215162500q^{10} $\\
    && $- 495925634096339/34426000q^9 + 661881614799897/43032500q^8 $\\
    && $-990040545191433/86065000q^7 + 532195382296631/107581250q^6 $\\
    && $-68592817731587/107581250q^5 - 18016218390139/86065000q^4 $\\
    && $ + 3138662105287/86065000q^3 - 179346320519/86065000q^2 $\\
    && $+ 100696210283/860650000q + 345858921/215162500)$\\ 
    $h_5(q)$&=& $-(4703529938733/3765343750q^{17} - 312187952900811/30122750000q^{16}$\\
    && $+126623998816272/1882671875q^{15} - 452217452448/1536875q^{14} $ \\
    && $ + 1607553611083881/1506137500q^{13} - 23225757911147049/7530687500q^{12}$\\
    && $+ 13577004998680164/1882671875q^{11} - 102322811958017807/7530687500q^{10} $\\
    && $+ 60575199001615299/3012275000q^9 - 138239622154714851/6024550000q^8$\\
    && $+ 144057800854658571/7530687500q^7 - 77960187940027919/7530687500q^6$\\
    && $+ 43589751452773401/15061375000q^5 - 162453007276159/1506137500q^4$\\
    && $- 179455823157787/1506137500q^3 + 169143999630681/7530687500q^2$\\
    && $ -20121321798503/15061375000q + 352512529837/30122750000)$\\
    $i_5(q)$&=& $-(-1567843312911/3765343750q^{17} + 26653336319487/7530687500q^{16}$\\
    && $- 86966693521859/3765343750q^{15} + 3148708883421/30737500q^{14}$ \\
    && $ - 563932915933193/1506137500q^{13} + 4122355705949559/3765343750q^{12}$ \\
    && $- 19526227867469373/7530687500q^{11} + 37360695285761837/7530687500q^{10} $\\
    && $- 16930208295615617/2259206250q^9 + 3303121075154424/376534375q^8$\\
    && $- 57158066913069257/7530687500q^7 + 32953970097517957/7530687500q^6$\\
    && $- 2598615541113562/1882671875q^5 + 6072469413529/60245500q^4$\\
    && $+ 246544595727319/4518412500q^3 - 21335413566293/1882671875q^2$\\
    && $ +5106134343567/7530687500q - 24469239909/3765343750)$\\
    $j_5(q)$&=& $-(3601089204639/15061375000q^{17} - 61218516478863/30122750000q^{16}$\\
    && $+ 49926764519649/3765343750q^{15} - 451845801429/7684375q^{14}$ \\
    && $ + 323616403589463/1506137500q^{13} - 730079188013933/7530687500q^{12} $\\
    && $+ 2799402050398322/1882671875q^{11} - 5353236974494508/1882671875q^{10}$\\
    && $+ 6463254746914883/1506137500q^9 - 30221823807366039/6024550000q^8$\\
    \end{tabular}

    \begin{tabular}{lcl}
    && $+ 32609306647168767/7530687500q^7 - 9353927650474821/3765343750q^6$\\
    && $+ 11642449338846181/15061375000q^5 - 15305412398327/301227500q^4$\\
    && $- 47792154921933/1506137500q^3 + 11564875919008/1882671875q^2$\\
    && $ -1292955862101/3765343750q + 127740623617/30122750000)$\\
    $k_5(q)$&=& $-(-3601089204639/30122750000q^{17} + 14666955541963/15061375000q^{16}$\\
    && $- 23688035104319/3765343750q^{15} + 334691467723/12295000q^{14}$ \\
    &&$  - 295534691350897/3012275000q^{13} + 2113643539915249/7530687500q^{12}$ \\
    && $- 9774053665697467/15061375000q^{11} + 6809946224333299/5648015625q^{10}$\\
    && $- 31633894300642129/18073650000q^9 + 2931470082058359/1506137500q^8$\\
    && $- 23470506685652367/15061375000q^7 + 17610921775161487/22592062500q^6$\\
    && $-672930270107769/3765343750q^5 - 63401012763611/9036825000q^4$\\
    && $+ 25429230732089/3012275000q^3 - 17299182621087/15061375000q^2$\\
    && $+ 4999986082313/90368250000q - 14793736609/22592062500)$
    \end{tabular}}

and $q$ such that $309q^{18} - 2781q^{17} + 18399q^{16} - 84156q^{15} + 313080q^{14} - 939456q^{13} + 2287664q^{12}- 4522686q^{11} + 7125279q^{10} - 8783755q^9 + 8166539q^8 - 5283806q^7 + 2051734q^6 - 317066q^5- 46820q^4 + 19554q^3 - 2121q^2 + 89q - 1=0$.

\subsection{}\label{5.12}$\omega^6_{5,6,i}=[f_6(q)x^3y^3+g_6(q)x^2y^4+h_6(q)xy^5 +i_6(q)y^6]dx + [j_6(q)x^6 + k_6(q)x^5y + l_6(q)x^4y^2 +x^3y^3]dy$ 

{\small    \begin{tabular}{lcl}
    $f_6(q)$&=& $-(-5/8pq^5 + 7/4pq^4 + 1/4q^5 + 13/8pq^3 - 5/8q^4 - 21/4pq^2 - q^3 + 3/8pq + 17/8q^2 + 3/2p$\\
    && $ + q - 5/8)$\\
    $g_6(q)$&=& $-(33/16pq^5 - 87/16pq^4 - 15/16q^5 - 99/16pq^3 + 9/4q^4 + 261/16pq^2 + 57/16q^3 + 33/16pq$\\
    && $- 27/4q^2 - 87/16p - 63/16q + 3/2)$\\
    $h_6(q)$&=& $-(-21/10pq^5 + 429/80pq^4 + 81/80q^5 + 267/40pq^3 - 189/80q^4 - 1281/80pq^2 - 303/80q^3$ \\
    &&$- 129/40pq + 537/80q^2 + 81/16p + 291/80q - 39/80)$\\
    $i_6(q)$&=& $-(7/10pq^5 - 7/4pq^4 - 7/20q^5 - 23/10pq^3 + 4/5q^4 + 26/5pq^2 + 13/10q^3 + 5/4pq - 11/5q^2$ \\
    &&$- 31/20p - 11/10q)$\\
    $j_6(q)$&=& $-(3/80pq^5 - 9/80pq^4 - 1/80q^5 - 9/80pq^3 + 1/20q^4 + 31/80pq^2 + 3/80q^3 + 23/80pq - 3/20q^2$ \\
    &&$- 1/16p - 1/80q + 1/20)$\\
    $k_6(q)$&=& $-(-3/40pq^5 + 3/16pq^4 + 3/80q^5 + 3/10pq^3 - 9/80q^4 - 51/80pq^2 - 9/80q^3 - 3/4pq + 21/80q^2 $\\
    &&$- 21/80p - 27/80q - 3/16)$\\
    $l_6(q)$&=&$ -(3/4p + 3/4q + 3/4)$
    \end{tabular}}

and $q$ such that for each $i=1,\cdots,6$: 

$i= \left \{ \begin{matrix}
        1, && q^2-q-1=0, && p =q-1 \\
        2, && q^2-q-1=0, && p =q+1 \\
        3, && q^2+q-1=0, && p =-q+1 \\
        4, && q^2+q-1=0, && p =q+1 \\
        5, && q^2-3q+1=0, && p =q-1 \\
        6, && q^2-3q+1=0, && p =-q+1 
    \end{matrix} \right .$

\subsection{}\label{5.13} $\omega^6_{6,1,i}=[f_1(r)x^4y^2+g_1(r)x^3y^3+h_1(r)x^2y^4+i_1(r)xy^5 +j_1(r)y^6]dx + x^6dy$

{\small    \begin{tabular}{lcl}
        $f_1(r)$&=& $ -(111/11r^{11} - 612/11r^{10} + 1902/11r^9 - 3972/11r^8 + 6225/11r^7 - 7536/11r^6 + 7644/11r^5$\\
        &&$- 6480/11r^4 + 4377/11r^3 - 2124/11r^2 + 702/11r - 60/11)$\\
        $g_1(r)$&=& $-(-540/11r^{11} + 2976/11r^{10} - 9248/11r^9 + 19308/11r^8 - 30236/11r^7 + 36556/11r^6$ \\
        &&$- 37012/11r^5 + 31320/11r^4 - 21128/11r^3 + 10244/11r^2 - 3384/11r + 448/11)$\\
        $h_1(r)$&=& $-(1065/11r^{11} - 5865/11r^{10} + 18219/11r^9 - 38013/11r^8 + 59469/11r^7 - 71805/11r^6 $\\
        &&$+ 2603/11r^5 - 61377/11r^4 + 41385/11r^3 - 20061/11r^2 + 6627/11r - 996/11)$\\
        $i_1(r)$&=& $-(-954/11r^{11} + 5250/11r^{10} - 1482r^9 + 3090r^8 - 4830r^7 + 5826r^6 - 5886r^5 + 4974r^4$\\
        && $- 3354r^3 + 1626r^2 - 5910/11r + 954/11)$\\
        
        $j_1(r)$&=& $-(318/11r^{11} - 159r^{10} + 5429/11r^9 - 11313/11r^8 + 17672/11r^7 - 21301/11r^6 + 21511/11r^5$\\
        && $- 18177/11r^4 + 12260/11r^3 - 5945/11r^2 + 1965/11r - 335/11)$
    \end{tabular}}

and such that  for each $i=1,\cdots,6$:
   $$i= \left \{ \begin{matrix}
        1, && r^4-3r^3+4r^2-2r+1=0, && q=r^3-2r^2+2r\\
        2, && r^4-3r^3+4r^2-2r+1=0, && q=r^2-r+1\\
        3, && r^4-2r^3+4r^2-3r+1=0, && q=-r^3+2r^2-3r+2\\
        4, && r^4-2r^3+4r^2-3r+1=0, && q=-r^3+r^2-2r+1\\
        5, && r^4-r^3+r^2-r+1, && q=-r^2+r\\
        6, && r^4-r^3+r^2-r+1, && q=r^3-r^2+r
    \end{matrix} \right .$$

\subsection{}\label{5.14}$\omega^6_{6,2}=[f_2(r)x^3y^3+g_2(r)x^2y^4+h_2(r)xy^5 +i_2(r)y^6]dx + [j_2(r)x^6 +x^5y]dy$

{\small    \begin{tabular}{lcl}
        $f_2(r)$&=&  $-(-211753/509775r^{11} + 1487389/169925r^{10} - 39488398/509775r^9 + 383308229/1019550r^8$\\ 
        &&$- 847379353/1019550r^7 + 1319404277/1019550r^6 - 1417830503/1019550r^5$\\
        &&$+ 441870419/339850r^4 - 867839501/1019550r^3 + 66039984/169925r^2 $\\
        &&$- 13645281/169925r + 4630952/509775)$\\
        $g_2(r)$&=&$ -(3369/1942r^{11} - 1767537/48550r^{10} + 7777346/24275r^9 - 149818851/97100r^8 $\\
        &&$+ 324532721/97100r^7 - 498940813/97100r^6 + 262863707/48550r^5 - 122298837/24275r^4$ \\
        &&$+ 310383701/97100r^3 - 68687683/48550r^2 + 6108456/24275r - 2086429/971001)$\\
        $h_2(r)$&=& $-(-1866449/849625r^{11} + 39088648/849625r^{10} - 343021373/849625r^9 $\\
        &&$+3290884181/1699250r^8 - 3526864102/849625r^7 + 5385766637/849625r^6$ \\
        &&$ - 5614394913/849625r^5 + 10420864211/1699250r^4 - 3247193192/849625r^3$\\
        &&$+ 1410585972/849625r^2 - 225894707/849625r + 35185247/1699250)$
        \end{tabular}
    
    \begin{tabular}{lcl}
        $i_2(r)$&=& $-(21406/24275r^{11} - 2686471/145650r^{10} + 23534869/145650r^9 - 56311732/72825r^8$\\
        &&$+ 479664209/291300r^7 - 729395167/291300r^6 + 251735073/97100r^5 - 174919307/72825r^4$\\
        &&$+ 430979153/291300r^3 - 46219883/72825r^2 + 6842413/72825r - 141879/19420)$\\
        $j_2(r)$&=&$ -(-22661/5097750r^{11} + 232156/2548875r^{10} - 3950857/5097750r^9 + 12078153/3398500r^8$\\
        &&$- 11336711/1699250r^7 + 7864483/849625*r^6 - 82517149/10195500r^5 + 18229171/2548875r^4$\\
        &&$- 14922823/5097750r^3 + 2883173/5097750r^2 + 1491001/2548875r + 13249/5097750)$
    \end{tabular}}

and $r$ such that $2r^{12} - 42r^{11} + 370r^{10} - 1785r^9 + 3887r^8 - 6020r^7 + 6413r^6 - 6020r^5 + 3887r^4 - 1785r^3 + 370r^2 - 42r + 2=0$.

\subsection{}\label{5.15}$\omega^6_{6,3}=[f_3(r)x^3y^3+g_3(r)x^2y^4+h_3(r)xy^5 +i_3(r)y^6]dx + [j_3(r)x^6 +k_3(r)x^5y +x^4y^2]dy$

{\small    \begin{tabular}{lcl}
        $f_3(r)$&=& $ -(-178679/1309750r^{11} + 99656/130975r^{10} + 397257/1309750r^9 - 366708/50375r^8 $ \\
        &&$- 65040184/654875r^7 + 8373093/21125r^6 - 1314001349/1964625r^5 + 2576672053/3929250r^4$\\
        &&$- 10039063/30225r^3 + 176894441/3929250r^2 + 19323407/1309750r - 1210419/261950)$\\
        $g_3(r)$&=&$ -(280224/654875r^{11} - 12412809/5239000r^{10} - 5462343/5239000r^9 + 229446/10075*r^8$ \\
        &&$+ 818323767/2619500r^7 - 1611557229/1309750r^6 + 6934773/3380r^5 - 5189950089/2619500r^4$\\
        &&$+ 390315957/403000r^3 - 539932383/5239000r^2 - 14370627/261950r + 2534041/209560)$\\
        
        \end{tabular}

        \begin{tabular}{lcl}
        $h_3(r)$&=&$ -(-2826651/6548750r^{11} + 15588057/6548750r^{10} + 3604563/3274375r^9 - 11551473/503750r^8$\\
        &&$- 206675652/654875r^7 + 8080695543/6548750r^6 - 13393623791/6548750r^5$  \\
        &&$+ 12835828463/6548750r^4 - 237065263/251875r^3 + 567257641/6548750r^2$ \\
        &&$+ 192527703/3274375r - 66420573/6548750)$\\
        $i_3(r)$&=& $-(942217/6548750r^{11} - 10364387/13097500r^{10} - 4944339/13097500*r^9 + 147903/19375r^8 $\\
        &&$+ 689594281/6548750r^7 - 268314767/654875r^6 + 2658380611/3929250r^5$ \\
        &&$- 6339239356/9823125r^4 + 926355493/3022500r^3 - 1003022297/39292500r^2$ \\
        &&$- 4268411/211250r + 114561287/39292500)$\\
        $j_3(r)$&=& $-(24411/6548750r^{11} - 564693/26195000r^{10} - 142239/26195000r^9 + 52119/251875r^8$ \\
        &&$+ 7027959/2619500r^7 - 74996953/6548750r^6 + 254071327/13097500r^5$  \\
        &&$- 246405651/13097500r^4 + 17960849/2015000r^3 - 7678599/26195000r^2$\\
        &&$- 2972518/3274375r + 3385091/78585000)$\\
        $k_3(r)$&=&$ -(-24411/3274375r^{11} + 268521/6548750r^{10} + 140247/6548750r^9 - 15651/38750r^8$\\
        &&$- 17907618/3274375r^7 + 27915279/1309750r^6 - 44032883/1309750r^5 + 96977867/3274375r^4$\\
        &&$- 2664749/251875r^3 - 9038099/3274375r^2 + 511581/211250r + 2332689/3274375)$\\
    \end{tabular}}

and $r$ such that $3r^{12} - 18r^{11} + 165r^9 + 2118r^8 - 9660r^7 + 18092r^6 - 19575r^5 + 11798r^4 - 2670r^3 - 440r^2 + 187r - 17=0$.

\subsection{}\label{5.16}$\omega^6_{6,4,i}=[f_4(r)x^2y^4+g_4(r)xy^5 +h_4(r)y^6]dx + [i_4(r)x^6 +j_4(r)x^5y +x^4y^2]dy$

{\small    \begin{tabular}{lcl}
        $f_4(r)$&=& $ -(19/75qr^7 - 467/150qr^6 - 13/240r^7 + 5/3q^3r^3 + 3643/150qr^5 + 311/600r^6 - 131/10q^3r^2$\\
        &&$- 29/25q^2r^3 - 6057/50qr^4 - 3313/1200r^5 + 59/6q^3r + 683/75q^2r^2 + 17279/150qr^3$ \\
        &&$+ 19577/1200r^4 - 359/30q^3 - 169/25q^2r - 20021/150qr^2 - 10421/400r^3 + 623/75q^2$ \\
        &&$+ 4919/150qr + 413/16r^2 - 1111/75q - 2041/120r + 77/16)$\\
        $g_4(r)$&=& $-(-111/250qr^7 + 682/125qr^6 + 19/200r^7 - 1097/375q^3r^3 - 10641/250qr^5 - 1363/1500r^6$\\
        &&$+ 8623/375q^3r^2 + 152/75q^2r^3 + 159251/750qr^4 + 14509/3000r^5 - 2159/125q^3r - 398/25q^2r^2$ \\
        &&$- 151453/750qr^3 - 85781/3000r^4 + 2626/125q^3 + 892/75q^2r + 58499/250qr^2 + 137083/3000r^3$\\
        &&$- 1094/75q^2 - 21574/375qr - 135811/3000r^2 + 6501/250q + 44777/1500r - 5577/1000)$\\
        $h_4(r)$&=&$ -(1831/9000qr^7 - 11249/4500qr^6 - 49/1125r^7 + 302/225*q^3r^3 + 175523/9000qr^5$  \\
        &&$+ 937/225r^6 - 2374/225q^3r^2 - 1043/1125q^2r^3 - 58381/600qr^4 - 2492/1125r^5 + 1784/225q^3r $ \\
        &&$+ 8197/1125q^2r^2 + 33317/360qr^3 + 2948/225r^4 - 241/25q^3 - 2051/375q^2r - 965221/9000qr^2$  \\
        &&$- 4714/225r^3 + 7532/1125q^2 + 118729/4500qr + 934/45r^2 - 107389/9000q$\\
        &&$- 15407/1125r + 3997/2250)$\\
        $i_4(r)$&=& $-(11/9000qr^7 - 67/4500qr^6 - 1/3600r^7 + 11/1125q^3r^3 + 1051/9000qr^5 + 23/9000r^6$\\
        &&$- 173/2250q^3r^2 - 2/225q^2r^3 - 1769/3000qr^4 - 229/18000r^5 + 127/2250q^3r + 16/225q^2r^2$ \\
        &&$+ 5101/9000qr^3 + 1421/18000r^4 - 17/250q^3 - 4/75q^2r - 6029/9000qr^2 - 2443/18000r^3$\\
        &&$ + 11/225q^2 + 803/4500qr + 2411/18000r^2 - 761/9000q - 857/9000r + 211/18000)$\\
        $j_4(r)$&=& $-(-7/500qr^7 + 43/250qr^6 + 3/1000r^7 - 7/75q^3r^3 - 671/500qr^5 - 43/1500r^6 + 11/15q^3r^2$\\
        &&$+ 26/375q^2r^3 + 2009/300qr^4 + 457/3000r^5 - 41/75q^3r - 68/125q^2r^2 - 637/100qr^3$ \\
        &&$- 541/600r^4 + 2/3q^3 + 146/375q^2r + 3697/500qr^2 + 289/200r^3 - 58/125q^2 - 1369/750qr$ \\
        &&$- 859/600r^2 + 413/500q + 1421/1500r - 71/3000)$
    \end{tabular}}

and $r$ is a parameter such that for each $i=1,2,3$:

$i= \left \{ \begin{matrix}
        1, && r^4 - 5r^3 + 52r^2 - 5r + 1=0, && 8q^2+q(-r^3 + 6r^2 - 58r - 1) + 8r=0 \\
        2, && r^4 - 8r^3 + 7r^2 - 8r + 1=0, && q^2+q(r^2 -r + 1) + r^2=0  \\
        3, && r^4 - 8r^3 + 7r^2 - 8r + 1=0, && q^2 +q(-r^3 + 8r^2  - 6r) + 1=0
    \end{matrix} \right .$

\subsection{}\label{5.17}$\omega^6_{6,5,i}=[f_5(r)x^3y^3+g_5(r)x^2y^4+h_5(r)xy^5 +i_5(r)y^6]dx + [j_5(r)x^6 +k_5(r)x^5y +l_5(r)x^4y^2$

{\small \begin{tabular}{lcl}
   $f_5(r)$&=&  $-(-1610951088177/865133593750r^{23} + 18378123031983/865133593750r^{22}$\\
        &&$+ 452195080667337/865133593750r^{21} - 10804121077881663/1730267187500r^{20}$\\
        &&$+ 76624272604592/432566796875r^{19} + 1024391296988695541/5190801562500r^{18}$\\
        &&$- 233596065819442979/432566796875r^{17} - 7218590085464114213/5190801562500r^{16}$\\
        &&$+ 4163187265419527173/519080156250r^{15} - 6411390285114593331/865133593750r^{14}$ \\
        &&$- 29610681684913996051/1297700390625r^{13} + 360999388151545693667/5190801562500r^{12}$\\
        &&$- 68814556528088246093/865133593750r^{11} + 99357115724477627/2525937500r^{10}$\\
        &&$+ 31571337711573539/21994921875r^9 - 52629664391798657609/5190801562500r^8$\\
        &&$+ 1370395419083864501/432566796875r^7 + 1199853768803905921/2595400781250r^6$\\
        &&$- 178736631843091009/519080156250r^5 + 29469456231323593/1730267187500r^4$\\
        &&$+ 26442382532764769/2595400781250r^3 - 1695886351115797/1730267187500r^2 $\\
        &&$- 6531799231439/173026718750r + 4554201340293/1730267187500)$\\
        $g_5(r)$&=& $ -(11319536172861/3460534375000r^{23} - 256046767416999/6921068750000r^{22}$\\
        &&$- 127587057096897/138421375000r^{21} + 75282668460305739/6921068750000r^{20}$\\
        &&$+ 1267203445089941/1730267187500r^{19} - 2396865876716053793/6921068750000r^{18}$\\
        &&$+ 791379497734313669/865133593750r^{17} + 4366210897208720373/1730267187500r^{16}$\\
        &&$- 5978228998221359326/432566796875r^{15} + 81241339660498398209/6921068750000r^{14}$\\
        &&$+ 17690670107823993087/432566796875r^{13} - 816590609337779966993/6921068750000r^{12}$\\
        &&$+ 55938956255672203317/432566796875r^{11} - 408720077347552201061/6921068750000r^{10}$\\
        &&$- 20859408238109159973/3460534375000r^9 + 7217944489498948457/432566796875r^8$ \\
        &&$- 15192400206981458849/3460534375000r^7 - 6573007730691797763/6921068750000r^6$\\
        &&$+ 886105494850884811/1730267187500r^5 - 81756081570093761/6921068750000r^4$ \\
        &&$-2158398889517819/138421375000r^3 + 8785966533190851/6921068750000r^2$\\
        &&$ + 76390607274063/1730267187500r - 1267590361461/346053437500)$\\
        $h_5(r)$&=& $-(-11319536172861/8651335937500r^{23} + 249215412084021/17302671875000r^{22}$\\
        &&$+ 806812220438451/2162833984375r^{21} - 3666727761813819/865133593750r^{20}$\\
        &&$- 6795004455962471/4325667968750r^{19} + 119426364887983401/865133593750r^{18}$\\
        &&$- 2805414023374674377/8651335937500r^{17} - 19152480864879481339/17302671875000r^{16} $\\
        &&$+ 11234274274914629283/2162833984375r^{15} - 54196862877293099339/17302671875000r^{14}$\\
        &&$- 7483157406988457267/432566796875r^{13} + 181660644832219853589/4325667968750r^{12}$\\
        &&$- 338281852378687236117/8651335937500r^{11} + 51096857544938054291/4325667968750r^{10}$\\
        &&$+ 223214113695753233/36658203125r^9 - 17147036813286646431/3460534375000r^8$\\
        &&$+ 625636245104905847/2162833984375r^7 + 8471508690345226621/17302671875000r^6$\\
        &&$- 595285484368391703/8651335937500r^5 - 74469384003231421/4325667968750r^4$\\
        &&$+ 9369559759681113/4325667968750r^3 + 176230569892527/865133593750r^2 $\\
        &&$+ 19599446489301/2162833984375r + 3591378323229/17302671875000)$\\
        $i_5(r)$&=& $-(1855653315297/8651335937500r^{22} - 20412186468267/8651335937500r^{21}$\\
        &&$- 2116768542747399/34605343750000r^{20} + 2402539153303137/3460534375000r^{19}$\\
        &&$+ 27220618527791279/103816031250000r^{18} - 195668609769111237/8651335937500r^{17}$\\
        &&$+ 5500561885957077421/103816031250000r^{16} + 4711519766687151919/25954007812500r^{15}$\\
        &&$- 367398013847629564/432566796875r^{14} + 26461541987468639887/51908015625000r^{13}$\\
        &&$+ 293597001766476994211/103816031250000r^{12} - 29669475006661122991/4325667968750r^{11}$\\
        &&$+ 221480219931687309437/34605343750000r^{10} - 5118195084913636121/2595400781250r^9$ \\
        &&$- 98822193808376163911/103816031250000r^8 + 13835645356715673011/17302671875000r^7$\\ 
        &&$- 1475938109917667311/25954007812500r^6 - 975202048472531843/12977003906250r^5$ \\
        &&$+ 82550717774880209/6921068750000r^4 + 59142407372359501/25954007812500r^3$ \\
        &&$- 13660511178748951/34605343750000r^2 - 50615608440419/3460534375000r$\\
        &&$+ 26915307014971/34605343750000)$\\
        \end{tabular}
    
    \begin{tabular}{lcl}
        $j_5(r)$&=& $-(-1739586755043/17302671875000r^{23} + 40765512355533/34605343750000r^{22}$\\
        &&$+ 482501946309327/17302671875000r^{21} - 479010283366491/1384213750000r^{20}$\\
        &&$+ 234009377632568/2162833984375r^{19} + 74129973854836557/6921068750000r^{18}$\\
        &&$- 140146688850695929/4325667968750r^{17} - 1177077918699471481/17302671875000r^{16}$\\
        &&$+ 994909245678624372/2162833984375r^{15} - 17883717897271396627/34605343750000r^{14}$\\
        &&$- 2031868608769161197/1730267187500r^{13} + 143828927888251927473/34605343750000r^{12}$\\
        &&$- 11277798913525965392/2162833984375r^{11} + 102157313314039825677/34605343750000r^{10}$\\
        &&$- 2874359753597865541/17302671875000r^9 - 2321272670930391339/3460534375000r^8$\\
        &&$+ 4683193009175091929/17302671875000r^7 + 618136224417732193/34605343750000r^6$\\
        &&$- 239578662684526177/8651335937500r^5 + 83946642336032313/34605343750000r^4$\\
        &&$+ 13964883302673613/17302671875000r^3 - 712399192452759/6921068750000r^2 $\\
        &&$- 6066413067936/2162833984375r + 5774407676991/17302671875000)$\\
        $k_5(r)$&=& $-(1653829643799/1730267187500r^{23} - 190945426026429/17302671875000r^{22}$\\
        &&$- 1154417688679869/4325667968750r^{21} + 28055092906275201/8651335937500r^{20}$\\
        &&$- 17648761283251/34605343750r^{19} - 878306311963967517/8651335937500r^{18}$\\
        &&$+ 2517338884499729437/8651335937500r^{17} + 11826005492692577309/17302671875000r^{16}$\\
        &&$- 18304330908033080589/4325667968750r^{15} + 14897667929994298559/3460534375000r^{14}$\\
        &&$+ 24798132351724088087/2162833984375r^{13} - 323895804120521605503/8651335937500r^{12}$\\
        &&$+ 387433719064862702607/8651335937500r^{11} - 204519363468248907833/8651335937500r^{10}$\\
        &&$+ 7883694273491317/86513359375r^9 + 102129056491673792181/17302671875000r^8$\\
        &&$- 4514590001928719992/2162833984375r^7 - 4005458641262114501/17302671875000r^6$ \\
        &&$+ 1974873160794456957/8651335937500r^5 - 5323158363264803/346053437500r^4$ \\
        &&$- 15137391915488103/2162833984375r^3 + 6892888476556899/8651335937500r^2$ \\
        &&$+ 21488468917959/865133593750r - 45463760972523/17302671875000$\\
        $l_5(r)$&=& $-(-1653829643799/1730267187500r^{23} + 38038081807377/3460534375000r^{22}$\\
        &&$+ 231298797080223/865133593750r^{21} - 22358031253906761/6921068750000r^{20}$\\
        &&$+ 1276524464563213/3460534375000r^{19} + 140460052555015447/1384213750000r^{18}$\\
        &&$- 99113461384353163/346053437500r^{17} - 961067923518396833/1384213750000r^{16}$\\
        &&$+ 7259909337113165911/1730267187500r^{15} - 14305848143818328297/3460534375000r^{14}$\\
        &&$- 8015132224926351659/692106875000r^{13} + 51075756867700475131/1384213750000r^{12}$\\
        &&$- 7506705617026219287/173026718750r^{11} + 2603637157898920601/117306250000r^{10}$\\
        &&$+ 952278199136959949/1730267187500r^9 - 8128216670119065781/1384213750000r^8$ \\
        &&$+ 1349559038410098497/692106875000r^7 + 185693391200930197/692106875000r^6$\\
        &&$- 97072896248137921/432566796875r^5 + 90152392073179051/6921068750000r^4$ \\
        &&$+ 6187252618611467/865133593750r^3 - 5482075040910849/6921068750000r^2 $\\
        &&$- 79136350624947/3460534375000r + 924487968159/276842750000)$
    \end{tabular}}

and $r$ such that for each $i=1,\cdots,3$: 

$i= \left \{ \begin{matrix}
        1, && r^8 + 14r^7 - 52r^6 - 168r^5 + 419r^4 - 168r^3 - 52r^2 + 14r + 1=0\\
        2, && r^8 - 22r^7 + 74r^6 + 130r^5 - 641r^4 + 666r^3 - 180r^2 - 36r + 9=0\\
        3, && 9r^8 - 36r^7 - 180r^6 + 666r^5 - 641r^4 + 130r^3 + 74r^2 - 22r + 1=0
    \end{matrix} \right .$

\section{Appendix}

For the readers' better comprehension of the procedure employed in the preceding section, we present an example of a SageMath script. This script pertains to case $\mathcal{T}_\h=2\cdot R_1 + 1\cdot R_3 + 1 \cdot R_5$.

\bigskip

\footnotesize
\begin{verbatim}
R=PolynomialRing(QQ, 'a6,a5,a4,p')
R.inject_variable()

S=PolynomialRing(R,'z')
S.inject_variable()

A=a6*z^6+a5*z^5+a4*z^4
B=1
A,B

E=5*[0]
E[1] = (A + B)(z=1)
E[2] = (A+B).derivative(z)(z=1)
E[3] = (A + p*B)(z=p)
E[4] = (A+p*B).derivative(z)(z=p)
L=[E[i] for i in range(1,5)]
for l in L:
	print(l)
	print("---")
	
restrictions = [
p,(p-1)
]
J=R.ideal(L)
J

I=J
for i in range(0,len(restrictions)):
	I=I.saturation(R.ideal(restrictions[i]))
	print("restriction: {0}, saturation: {1}".format(restrictions[i],I[1]))
	I=I[0]
I.gens()

I.dimension()

I.vector_space_dimension()

dec=I.primary_decomposition(algorithm = 'gtz')

for c in dec:
	for p in c.gens():
		print(p)
		print("--")
	print("---------------------")

print(len(dec))
\end{verbatim}

\bigskip

\normalsize

Note that the lines "I.dimension()" and "I.vector\_ space\_ dimension()" are not necessary. However, they provides interesting information about the problem.
%%%%%%%%%%%%%%%%%%%%%%%%%%%%%%%%%%%%%%%%%%%%%%%%%%%%%%%%%%%%%%%%%%%%%%%%%%%%%%%%%%%%%%%%%%%%%%%%%%%%%%%%%%%%%%%%%%%%%%%%%%%%%%%%%%%%%%%%%%%%%%%%%%%%%%%%%%%%%%%%%%%%%%%%%%%%%%%%%%%%

\end{document}